\documentclass[a4paper,11pt,final]{article}

\usepackage{graphicx}
\usepackage{amssymb,amsmath,amsthm,mathrsfs,xfrac}
\usepackage[shortlabels]{enumitem}
\usepackage{a4wide}
\usepackage{hyperref}
\usepackage{array}
\usepackage{colonequals}
\usepackage{dsfont}
\usepackage{tikz}

\numberwithin{equation}{section}
\allowdisplaybreaks[4]

\newtheorem{proposition}{Proposition}[section]
\newtheorem{theorem}[proposition]{Theorem}

\newtheorem{lemma}[proposition]{Lemma}
\newtheorem{definition}[proposition]{Definition}
\newtheorem{corollary}[proposition]{Corollary}
\newtheorem{remark}[proposition]{Remark}

\newenvironment{proofof}[1]{\smallskip\noindent{\textbf{Proof~of~#1.}}%
  \hspace{1pt}}{\hspace{-5pt}{\nobreak\quad\nobreak\hfill\nobreak%
    $\square$\vspace{2pt}\par}\smallskip\goodbreak}

\newcommand{\Id}{\mathinner{\mathrm{Id}}}
\newcommand{\pint}[1]{\mathaccent23{#1}}
\newcommand{\C}[1]{\mathbf{C}^{#1}}
\newcommand{\AC}{\mathbf{AC}}
\newcommand{\ACR}{\mathbf{ACR}}
\newcommand{\Cc}[1]{\mathbf{C}_c^{#1}}

\newcommand{\BV}{\mathbf{BV}}

\newcommand{\Lloc}[1]{{\mathbf{L}_{\mathbf{loc}}^{#1}}}

\newcommand{\modulo}[1]{{\left|#1\right|}}
\newcommand{\norma}[1]{{\left\|#1\right\|}}
\newcommand{\caratt}[1]{\mathop{{\mathds{1}_{\strut#1}}}}
\newcommand{\reali}{{\mathbb{R}}}
\newcommand{\naturali}{{\mathbb{N}}}

\newcommand{\complessi}{{\mathbb{C}}}
\newcommand{\tv}{\mathinner{\rm TV}}
\newcommand{\spt}{\mathop{\rm spt}}
\newcommand{\sgn}{\mathop{\rm sgn}}

\newcommand{\wsto}{\stackrel{*}{\rightharpoonup}}

\renewcommand{\L}[1]{{\mathbf{L}^#1}}

\renewcommand{\epsilon}{\varepsilon}
\renewcommand{\phi}{\varphi}
\renewcommand{\theta}{\vartheta}

\renewcommand{\d}[1]{\mathinner{\mathrm{d}{#1}}}
\newcommand{\grad}{\mathinner{\mathop{{\rm grad}}}}
\renewcommand{\div}{\mathinner{\mathop{{\rm div}}}}

\newcommand{\Lin}{\mathbf{Lin}}
\newcommand{\LinC}{\mathbf{Lin}_{\mathbf{c}}}

\newcommand{\essinf}{\mathop{{\rm ess\,inf}}\limits}
\newcommand{\esssup}{\mathop{{\rm ess\,sup}}\limits}

\newcommand\rest[2]{{#1}_{\mkern 1mu \vrule height 2ex \mkern 2mu {#2}}}

\DeclareFontFamily{U}{mathx}{}
\DeclareFontShape{U}{mathx}{m}{n}{<-> mathx10}{}
\DeclareSymbolFont{mathx}{U}{mathx}{m}{n}
\DeclareMathAccent{\widecheck}{0}{mathx}{"71}

\title{Entropy Production in General Balance Laws}

\author{Rinaldo M.~Colombo$^1$ \and Vincent Perrollaz$^2$}

\begin{document}

\maketitle

\footnotetext[1]{INdAM Unit, University of Brescia,
  Italy. \texttt{rinaldo.colombo@unibs.it}}

\footnotetext[2]{Institut Denis Poisson, Universit\'e de Tours, CNRS
  UMR 7013, Universit\'e d'Orl\'eans, France\\
  \texttt{vincent.perrollaz@univ-tours.fr}}

\begin{abstract}

  \noindent Given a general scalar balance law, i.e., in several space
  dimensions and with flux and source both space and time dependent,
  we focus on the functional properties of the entropy production. We
  apply this operator to entropy solutions, to distributional
  solutions or to merely $\L\infty$ functions. Proving its analytical
  properties naturally leads to the projective tensor product of $\C1$
  spaces and to further natural extensions to space and time dependent
  or complex valued \emph{``entropies''}.  Besides various qualitative
  properties, this extended framework allows to obtain new
  quantitative {formul\ae} --- also by means of Fourier transforms ---
  that provide different representations of the entropy
  production. Remarkably, this operator also furnishes a
  \emph{faithful representation} of any $\L\infty$ function.

  \medskip

  \noindent\textit{2000~Mathematics Subject Classification:} 35L65;
  35L67; 46A32 

  \medskip


\end{abstract}


\section{Introduction}
\label{sec:Intro}

In the theory of general scalar balance law, i.e., of
\begin{equation}
  \label{eq:1}
  \left\{
    \begin{array}{l}
      \partial_t u + \div_x f (t,x,u) = g (t,x,u)
      \\
      u (0,x) = u_o (x)
    \end{array}
  \right.
  \qquad (t,x) \in \reali_+ \times \reali^n
\end{equation}
entropy has played a key role. Rooted in thermodynamics, its extension
to the case of~\eqref{eq:1} has helped the statements of uniqueness
and stability criteria. This paper is focused on the entropy
production operator, on its functional properties and on its use when
applied to distributional (i.e., not necessarily entropic) solutions
or even to general $\L\infty$ functions.

Indeed, the entropy production turns out to be a \emph{faithful
  representation} of any $\L\infty$ function. Then, the specific role
of {Kru\v zkov} entropies naturally stems out and general {formul\ae}
relating general entropies to {Kru\v zkov} entropies can be proved.

Different functional environments are considered which lead to
constructions first with different smoothness requirements and,
second, specifically adapted to non homogeneous settings. In
particular, an extension to complex valued entropies allows to obtain
a general representation formula by means of Fourier transform.

In this Introduction we first consider the by far simpler problem
\begin{equation}
  \label{eq:52}
  \left\{
    \begin{array}{l}
      \partial_t u + \partial_x f (u) = 0
      \\
      u (0,x) = u_o (x)
    \end{array}
  \right.
  \qquad (t,x) \in \reali_+ \times \reali \,,
\end{equation}
namely, in one space dimension, with $f$ independent of $(t,x)$ and
with no source term. The general case~\eqref{eq:1} is the scope of the
paper and it is in this framework that all statements hereafter are
proved in later sections. In the setting of~\eqref{eq:52}, for any
$u \in \L\infty (\reali;\reali)$ --- not necessarily related
to~\eqref{eq:52} --- the entropy production is the distribution valued
map
\begin{displaymath}
  \begin{array}{ccccc}
    \mathcal{M}_u
    & \colon
    & \C1 (\reali; \reali)
    & \to
    & \mathcal{D}' (\reali_+ \times \reali;\reali)
    \\
    &
    & E
    & \mapsto
    & \mathcal{M}_u (E)
  \end{array}
\end{displaymath}
where for all $\phi \in \mathcal{D} (I \times \reali;\reali)$ and for
all $F$ such that $(E,F)$ is an entropy - entropy flux pair, i.e.,
$F' = E' \, f' $,
\begin{eqnarray*}
  \left(\mathcal{M}_u (E)\right) (\phi)
  & \colonequals
  & \int_{I} \int_{\reali}
    \left(
    E\left(u (t,x)\right) \, \partial_t \phi (t,x)
    +
    F\left(u (t,x)\right) \cdot \partial_x \phi (t,x)
    \right) \d{x} \d{t}
  \\
  &
  & +
    \int_{\reali} E\left(u_o (x)\right) \, \phi (0,x) \d{x} \,.
\end{eqnarray*}
The map
$(f,u,u_o,E,\phi) \mapsto \left(\mathcal{M}_u (E)\right) (\phi)$ turns
out to be well defined, in spite of the non uniqueness of $F$ for
given $E$ and $f$. Remarkably, the topologies that allow to prove its
continuity impose to consider tensor product spaces and their
inductive topologies, see~\eqref{eq:42} Moreover,
Theorem~\ref{thm:Identifynotx} shows that for fixed $f$, $u_o$ the
correspondence $u \mapsto \mathcal{M}_u$ is injective on
$\L\infty (\reali; \reali)$. i.e.,
\begin{equation}
  \label{eq:26}
  \mbox{if } u,v \in \L\infty (\reali; \reali) \qquad\quad
  \mathcal{M}_u = \mathcal{M}_v \implies u = v \,.
\end{equation}
It is in this sense
that we state that $\mathcal{M}_u$ provides a \emph{faithful
  representation} of $u$.

Formally, {Kru\v zkov} entropies are obtained setting
$E_k (u) \colonequals \modulo{u-k}$, for $k \in \reali$, leading to
$\mu_k = \mathcal{M} (E_k)$ or, more precisely,
\begin{eqnarray*}
  \mu_k (\phi)
  & \!\!\!\! \colonequals \!\!\!\!\!
  & \int_{I} \! \int_{\reali^n}
    \! \left[
    \modulo{u (t,x) - k} \, \partial_t \phi (t,x)
    +
    \sgn \left(u (t,x) -k\right) \,
    \left(f \left(u (t,x)\right) - f (k)\right) \,
    \partial_x \phi (t,x)
    \right] \! \d{x} \d{t}
  \\
  &  \!\!\! \!\!\!
  & + \int_{\reali^n} \modulo{u_o(x) - k} \, \phi(0,x) \d{x} \,.
\end{eqnarray*}
Then, for every distributional solution $u$ to~\eqref{eq:52} --- but
it holds also in the case~\eqref{eq:1}, see Definition~\ref{def:ws}
--- we here prove the formula
\begin{equation}
  \label{eq:53}
  \mathcal{M}_u (E)
  =
  \dfrac12 \,
  \int_{\reali} E'' (k) \; \mu_k \d{k}
\end{equation}
In the general setting of~\eqref{eq:1}, the  integral
in~\eqref{eq:53} is a Dunford-Pettis integral in
$\mathcal{D}' (\reali;\reali)$, see
Appendix~\ref{prop:integral}.

A first consequence of~\eqref{eq:53} is the positivity of
$\mathcal{M}_u (E)$ for convex $E$. Besides, the
injectivity~\eqref{eq:26} can be restated saying, in the case of
distributional solutions, that the map $u \mapsto [k \mapsto \mu_k]$
is injective.

Again at a formal level, all the above can be extended to the case of
complex valued entropies and entropy fluxes. This trick allows to use
the machinery of Fourier transforms which, in turn, leads to a further
representation formula, see Proposition~\ref{prop:Fourier}. More than
that, we also get that if $E$ is $\C2$, $E''$ is a
tempered distribution and $\spt \widehat {E''}$ admits a cluster
point, then the translates of $E$ characterize distributional
solutions. Indeed, for $k \in \reali$, call
$\left(\mathcal{T}_k E\right) (u) = E (u-k)$ the $k$-translate of
$E$. If $u$, $v$ are distributional solutions to~\eqref{eq:52}, we
have the implications
\begin{displaymath}
  \forall \,k \in \reali \qquad
  \mathcal{M}_u \left(\mathcal{T}_kE \right) =
  \mathcal{M}_v \left(\mathcal{T}_kE\right)
  \quad \implies \quad
  \mathcal{M}_u = \mathcal{M}_v
  \quad \implies \quad
  u = v \,.
\end{displaymath}
Note that no convexity assumption on $E$ is required. In the spirit
of~\cite{MR4690554, MR1296003}, but in the present more general setting, the
above implications provide an initial description of \emph{``which''}
and \emph{``how many''} entropies single out a distributional
solution.

\goodbreak

The extension from~\eqref{eq:52} to~\eqref{eq:1}, in particular the
introduction of the $x$ dependence in $f$, needs particular care not
only from the technical point of view. The appearance of a $u$
derivative of $E$ in the expression of the entropy production,
see~\eqref{eq:21}, leads to a possible loss of regularity that menaces
the very meaningfulness of the various terms.

The first difficulty, overcame in \S~\ref{subsec:case-class-entr},
amounts to choose a functional framework where the equality
$\mu_k = \mathcal{M}_u (E_k)$ can be rigorously justified and the
representation formula~\eqref{eq:52} has now to be rewritten as
\begin{equation}
  \label{eq:54}
  \mathcal{M}_u (E)
  =
  \dfrac12 \, \int_{\reali} \mu_k \, \d{E' (k)} \,.
\end{equation}
Above, we have a Kurzweil-Stieltjes integral in
$\mathcal{D}' (\reali;\reali)$, see
Appendix~\ref{sec:integr-texorpdfstr}. Note that in~\eqref{eq:54}, $E$
is not required to be $\C2$, but only to admit a regulated derivative,
i.e., $E \in \ACR (\reali;\reali)$ in the sense of~\eqref{eq:20}, so
that the case $E_k (u) = \modulo{u-k}$ is rigorously comprised.

A further difficulty stems out from the appearance of $x$ derivatives
of $f$ and $F$ in~\eqref{eq:40}. Indeed, if both $f$ and $E$ are of
class $\C1$, the definition of the entropy flux associated to an
entropy $\partial_u F (t,x,u) = E' (u) \; \partial_u f (t,x,u)$ leads
to $F$ possibly lacking $\C1$ regularity in $(t,x)$. Clearly, this
issue also poses major obstacles to any regularity of
$(E,f) \mapsto F$. The tensor product
$\C1 (\reali\times\reali^n; \reali) \otimes \C1 (\reali;\reali)$
provides the natural algebraic setting. We thus consider maps $\zeta$
that are limit (with respect to the projective topology) of linear
combinations of functions $\zeta_1 (t,x) \, \zeta_2 (u)$, see
Appendix~\ref{sec:tensor-product-c1} for the rigorous treatment. Note
that any such function $\zeta$ admits the mixed second derivatives
$\partial^2_{tu} \zeta$ and $\nabla_x \partial_u\zeta$ and Schwarz
Lemma applies. It then comes natural to allow also $E$ to depend on
$(t,x,u)$ and, once we choose a specific entropy flux, e.g.,
$F (t,x,u) = \int_0^u \partial_u E (t,x,w) \; \partial_u f (t, x,w)
\d{w}$, the mapping $(f,E) \mapsto F$ turns out to be bilinear and
continuous in the projective tensor product
$\C1 (\reali_+ \times \reali^n; \reali) \otimes \C1 (\reali; \reali)$.

A further justification of the appearance of the tensor product
$\C1 (\reali_+ \times \reali^n; \reali) \otimes \C1 (\reali; \reali)$
is worth mentioning. If $X,Y,Z$ are vector spaces, denote by
$\Lin (X,Y)$ the set of linear maps defined in $X$ valued in
$Y$. Recall the natural isomorphism
\begin{equation}
  \label{eq:16}
  \Lin\left(X;\Lin (Y;Z)\right)
  \sim
  \Lin (X \otimes Y; Z)
\end{equation}
see~\cite[Theorem~39.1]{MR2296978}. From the results in Section~\ref{sec:main-resul} we have
\begin{displaymath}
  \mathcal{M}_u \in \Lin_c
  \left(
    \C1 (\reali;\reali);
    \LinC\left(\Cc1 (I\times\reali^n); \reali\right)
  \right)
\end{displaymath}
where $\LinC (X,Y)$ is the space of linear and continuous map defined
in $X$ valued in $Y$, $X$ and $Y$ being now topological vector
spaces. Section~\ref{sec:time-space-dependent} can be seen as the
proof that
\begin{displaymath}
  \mathcal{M}_u \in \Lin_c
  \left(
    \C1 (\reali;\reali) \otimes \Cc1 (I\times\reali^n); \reali
  \right)
\end{displaymath}
where $\C1 (\reali;\reali) \otimes \Cc1 (I\times\reali^n)$ is equipped
with the projective topology. Remark that general constructions
yielding this result are not available. Indeed, a general extension of
the isomorphism~\eqref{eq:16} consistent with the topological
structure is unavailable, unless in the case of nuclear
spaces~\cite[Chapter~50]{MR2296978}, not applicable to the present
setting. A general construction that may yield a result close to the
present one might be based on Schwartz Kernel
Theorem~\cite[Theorem~5.2.1]{MR1065993}. However, it would require to
restrict to $\Cc\infty$ entropies and fluxes.

This point of view, developed in
Section~\ref{sec:time-space-dependent}, rises a natural connection to
three entirely different and quite relevant formulations of
conservation laws, namely the kinetic one~\cite{MR1201239, MR2064166},
the measure valued one~\cite{MR775191} and the entropy process
one~\cite{MR1235201}, see also~\cite{MR2172652}. Indeed, they all
amount to substitute nonlinear conservation (or balance) laws with
linear evolutions of different nonlinear quantities. The basis for
this link is Formula~\eqref{eq:17} in
Theorem~\ref{thm:formulaC2exProp} --- which in a sense extends
both~\eqref{eq:53} and~\eqref{eq:54} --- and yields a linear evolution
of the map $(t,x,k) \mapsto \modulo{u (t,x) - k}$, obviously related
to {Kru\v zkov} entropies. The possible consequences of these facts
deserve a deeper analysis which we leave to forthcoming works.

Let us also note that the present developments allows a rephrasing of
known results, leaving open the question of providing alternative
proofs obtained \emph{within} the present framework. Here we refer in
particular to Corollary~\ref{cor:Panov} and
Corollary~\ref{cor:Kruzkov}.

Entropy production, when positive, has been well known to single out
the unique \emph{good} solution at least since the classical results
by {Kru\v zkov}~\cite{Kruzkov1970}. More recently, $\L1$ contractive
semigroups of distributional though non entropy solution have been
shown to actually exist, often justified by applications to vehicular
traffic, see~\cite{MR2807133, MR2300671}. This is a further motivation
for the present study of the entropy production operator in a general
functional setting.

\medskip

The next section presents the classical case of entropies depending
only on the $u$ variable: initially in the $\C1$ setting
(\S~\ref{subsec:first-setting}), then in the \emph{regulated} case
(\S~\ref{subsec:case-class-entr}) and, third, a complex valued
framework is considered (\S~\ref{subsec:fourier-transform}). Then,
Section~\ref{sec:time-space-dependent} considers entropies depending
also on $(t,x)$. Proofs are deferred to
Section~\ref{sec:technical-details}, the first to be presented being
those related to Section~\ref{sec:time-space-dependent}. The two final
appendixes gather general results of use in this paper but not easily
available.

\section{The Standard Entropy Operator}
\label{sec:main-resul}

Throughout, we let $I = \reali_+$ or $I= [0, T \mathclose[$ for an
arbitrary $T>0$. As usual, $\mathcal{D} (\Omega; \reali)$ stands for
$\Cc\infty (\Omega;\reali)$, for any $\Omega \subseteq \reali^m$.

The following spaces and topologies are used throughout. We define them
here and defer their main properties to
Appendix~\ref{sec:tensor-product-c1}.
\begin{equation}
  \label{eq:42}
  \mathscr{F} \colonequals \left\{ \zeta \in \C1 (\overline{I} \times
    \reali^n \times \reali;\reali) \colon
    \begin{array}{rcl@{}}
      \partial^2_{tu}\zeta
      & \in
      & \C0(\overline{I} \times \reali^n\times\reali;\reali)
      \\
      \grad_x \partial_u \zeta
      & \in
      & \C0(\overline{I} \times \reali^n \times \reali;\reali^n)
    \end{array}
  \right\}
\end{equation}
as well as $\mathscr{F}^n$ for functions valued in $\reali^n$.
Introduce, for any compact set
$K \subset I \times \reali^n \times \reali$, the seminorms
\begin{equation}
  \label{eq:45}
  p_k (\zeta) \colonequals \norma{\zeta}_{\C1 (K;\reali)} +
  \norma{\partial^2_{tu}\zeta}_{\C0 (K;\reali)} +
  \norma{\grad_x \partial_u \zeta}_{\C0 (K;\reali^n)} \,.
\end{equation}
For functions $\zeta$ defined in $I \times \reali^n \times \reali$,
introduce their support in the $t$ and $x$ variables as
\begin{equation}
  \label{eq:50}
  {\spt}_{tx} \zeta
  \colonequals
  \overline{\left\{
      (t,x) \in I \times \reali^n
      \colon
      \exists\, u \in \reali \quad \zeta (t,x,u) \neq 0
    \right\}} \,,
\end{equation}
the closure being in the topology of $I \times \reali^n$.  The space
\begin{equation}
  \label{eq:41}
  \mathscr{F}_c
  \colonequals
  \left\{
    \zeta \in \mathscr{F}  \colon
    {\spt}_{tx} \zeta \mbox{ is compact}
  \right\}
\end{equation}
is naturally equipped with the inductive topology, see
Lemma~\ref{lem:treves}.

Concerning the variables of the functions $\zeta \in \mathscr{F}$ or
$\zeta \in \mathscr{F}_c$ we use the notation $\zeta (t,x,u)$.

\subsection{Elementary Setting -- \texorpdfstring{$E \in \C1$}{E in
    C1}}
\label{subsec:first-setting}

With reference to~\eqref{eq:1}, we assume that the flux $f$ is in
$\mathscr{F}^n$ from~\eqref{eq:42}, the source $g$ is in
$\C0 (\overline{I} \times \reali^n \times \reali; \reali)$ and the
initial datum $u_o$ is in $\L\infty (\reali; \reali)$.

Recall the usual definition of entropy - entropy flux pair.

\begin{definition}
  \label{def:eefp}
  A pair of function $(E,F)$ is an \emph{entropy - entropy flux pair
    with respect to $f$} if for all $(t,x) \in I \times \reali^n$ and
  for a.e. $u \in \reali$,
  \begin{eqnarray}
    \nonumber
    E
    & \in
    & \AC (\reali; \reali)\,;
    \\
    \nonumber
    u \mapsto F (t,x,u)
    & \in
    & \AC (\reali;\reali)\,;
    \\
    \label{eq:24}
    \partial_u F (t,x,u)
    & =
    & E' (u) \; \partial_u f (t,x,u) \,.
  \end{eqnarray}
\end{definition}

\noindent The map $(E,f) \mapsto F$, for a particular choice of $F$,
is studied in Proposition~\ref{prop:fEF}.

Throughout, if $u \in \L\infty (I \times \reali^n; \reali)$ and
$u_o \in \L\infty (\reali^n, \reali)$, we denote
\begin{eqnarray}
  \label{eq:14}
  J_u
  & \colonequals
  & \left[
    \essinf_{(t,x) \in I \times \reali} u \,,\,
    \esssup_{(t,x) \in I \times \reali} u \right] \,;
  \\
  \label{eq:5}
  J_{u,u_o}
  & \colonequals
  & \left[
    \min\left\{
    \essinf_{(t,x) \in I \times \reali} u, \,
    \essinf_{x \in \reali} u_o
    \right\} \,,\,
    \max\left\{
    \esssup_{(t,x) \in I \times \reali} u, \,
    \esssup_{x \in \reali} u_o
    \right\}
    \right] \,.
\end{eqnarray}

The following definition quantifies the defect in
equality~\eqref{eq:24} when $u$ is only a \emph{distributional}
solution, see also~\cite[\S~4.5]{MR3468916}. However, it can be
considered also when $u$ is an arbitrary bounded measurable function.

\begin{definition}
  \label{def:muMnotx}
  Let $f \in \mathscr{F}^n$,
  $g \in \C0 (\overline{I} \times \reali^n\times \reali;\reali)$,
  $u_o \in \L\infty (\reali; \reali)$ and
  $u \in \L\infty (I \times \reali^n; \reali)$. With reference
  to~\eqref{eq:1}, define the \emph{entropy production distribution}
  \begin{displaymath}
    \begin{array}{ccccc}
      \mathcal{M}_u
      & \colon
      & \C1 (\reali; \reali)
      & \to
      & \mathcal{D}' (I \times \reali^n;\reali)
      \\
      &
      & E
      & \mapsto
      & \mathcal{M}_u (E)
    \end{array}
  \end{displaymath}
  where for all $\phi \in \mathcal{D} (I \times \reali^n;\reali)$ and
  for $F \in \mathscr{F}$ with $(E,F)$ as in
  Definition~\ref{def:eefp},
  \begin{eqnarray}
    \nonumber
    \left(\mathcal{M}_u (E)\right) (\phi)
    & \colonequals
    & \int_{I} \int_{\reali^n}
      \left(
      E\left(u (t,x)\right) \, \partial_t \phi (t,x)
      +
      F\left(t, x, u (t,x)\right) \cdot \grad_x \phi (t,x)
      \right) \d{x} \d{t}
    \\
    \nonumber
    &
    & -
      \int_{I} \int_{\reali^n}
      E' \left(u (t,x)\right)
      \; (\div_x f)\left(t, x, u (t,x)\right)
      \; \phi (t,x)
      \d{x} \d{t}
    \\
    \label{eq:21}
    &
    & +
      \int_{I} \int_{\reali^n}
      (\div_x F)\left(t, x, u (t,x)\right)
      \; \phi (t,x)
      \d{x} \d{t}
    \\
    \nonumber
    &
    & +
      \int_{I} \int_{\reali^n}
      E' \left(u (t,x)\right) \;
      g \left(t,x,u (t,x)\right) \; \phi (t,x)
      \d{x} \d{t}
    \\
    \nonumber
    &
    & +
      \int_{\reali^n} E\left(u_o (x)\right) \, \phi (0,x) \d{x} \,.
  \end{eqnarray}
\end{definition}

The use of $\div_x F$ in~\eqref{eq:21} is legitimate thanks to the
choice $f \in \mathscr{F}^n$ with $\mathscr{F}$ as in~\eqref{eq:42}
which allows to apply Proposition~\ref{prop:fEF} that ensures the
existence of at least one entropy flux $F \in \mathscr{F}^n$.

Apparently, $\mathcal{M}_u$ depends on the particular choice of the
flux $F$. Indeed, for any fixed entropy $E$, the corresponding flux
$F$ is not uniquely determined by~\eqref{eq:24}. The next proposition
also shows that, given problem~\eqref{eq:1}, $\mathcal{M}_u (E)$ is
well defined only as a function of $E$.

\begin{proposition}
  \label{prop:1notx}
  Let $f \in \mathscr{F}^n$,
  $g \in \C0 (\overline{I} \times \reali^n\times \reali;\reali)$,
  $u_o \in \L\infty (\reali; \reali)$ and
  $u \in \L\infty (I \times \reali^n; \reali)$. Then,
  \begin{enumerate}[label=\bf($\mathcal{\mathcal{M}}$\arabic*)]
  \item \label{item:9notx} Let $F_1$ and $F_2$ be two entropy fluxes
    related to $E$ as in~\eqref{eq:24} and call $\mathcal{M}_u^1$ and
    $\mathcal{M}_u^2$ the corresponding functions resulting from
    Definition~\ref{def:muMnotx}. Then,
    $\mathcal{M}_u^1 = \mathcal{M}_u^2$.
  \item \label{item:33} If $u$ is a strong solution to~\eqref{eq:1}
    with initial datum $u_o$, then $\mathcal{M}_u \equiv 0$.
  \item \label{item:1notx} $\mathcal{M}_u$ is linear and if
    $E_{\strut| J_{u,u_o}} = 0$, then $\mathcal{M}_u (E) =0$.
  \item \label{item:4notx} For all $E \in \C1 (\reali; \reali)$,
    $\mathcal{M}_u (E)$ is a distribution of order at most $1$
  \end{enumerate}
\end{proposition}

\begin{remark}
  \label{rem:Ispt}
  \rm The latter condition in~\ref{item:1notx} can be written as
  follows, for arbitrary $E_1,E_2 \in \C1 (\reali;\reali)$:
  \begin{displaymath}
    \mbox{For all } w \in J_{u,u_o} \quad E_1 (w) = E_2 (w)
    \quad \implies \quad
    \mathcal{M}_u (E_1) = \mathcal{M}_u (E_2) \,.
  \end{displaymath}
  Alternatively, one can say that $\mathcal{M}_u$ is \emph{supported}
  in $\C1 (J_{u,u_o}; \reali)$.
\end{remark}

\begin{remark}
  \label{rem:confusion}
  \rm Due to~\ref{item:4notx}, in the sequel we view
  $\mathcal{M}_u (E)$ also as a linear and continuous functional
  acting on $\Cc1 (I\times\reali^n;\reali)$, i.e.,
  $\mathcal{M}_u (E) \in \left(\Cc1 (I\times\reali^n;\reali)\right)'$.
\end{remark}

\begin{theorem}
  \label{thm:continuitynotx}
  For $\nu \in \naturali$, let $f,f_\nu \in \mathscr{F}^n$;
  $g, g_\nu \in \C0 (\overline{I} \times \reali^n\times
  \reali;\reali)$; $u_o, u_o^\nu \in \L\infty (\reali; \reali)$ and
  $u_\infty, u^\nu \in \L\infty (I \times \reali^n; \reali)$;
  $E,E_\nu \in \C1 (\reali;\reali)$. Assume that
  \begin{eqnarray*}
    f_\nu
    & \underset{\nu\to+\infty}{\longrightarrow}
    & f \mbox{ in } \mathscr{F}^n \mbox{ with respect to~\eqref{eq:45}};
    \\
    g_\nu
    & \underset{\nu\to+\infty}{\longrightarrow}
    & g \mbox{ uniformly on compact subsets of }
      I \times \reali^n \times \reali;
    \\
    u_o^\nu
    & \underset{\nu\to+\infty}{\longrightarrow}
    & u_o \mbox{ pointwise a.e.~in } \reali^n \mbox{ and }
      u_o^\nu \mbox{ is bounded in } \Lloc\infty (\reali^n; \reali);
    \\
    u^\nu
    & \underset{\nu\to+\infty}{\longrightarrow}
    & u_\infty \mbox{ pointwise a.e.~in } \reali \times \reali^n \mbox{ and }
      u^\nu \mbox{ is bounded in } \Lloc\infty (I \times \reali^n; \reali);
    \\
    E_\nu
    & \underset{\nu\to+\infty}{\longrightarrow}
    & E \mbox{ in } \C1 (\reali;\reali).
  \end{eqnarray*}
  \noindent Call $\mathcal{M}_{u_\infty}$ the entropy production
  functional introduced in Definition~\ref{def:muMnotx} and
  $\mathcal{M}_\nu$ the entropy production functional generated by the
  Cauchy problem
  \begin{displaymath}
    \left\{
      \begin{array}{l}
        \partial_t u + \div_x f_\nu (t,x,u) = g_\nu (t,x,u)
        \\
        u (0,x) = u_o^\nu (x)
      \end{array}
    \right.
  \end{displaymath}
  and by the function $u_\nu$, according to
  Definition~\ref{def:muMnotx}. Then,
  \begin{displaymath}
    \mathcal{M}_\nu (E_\nu)
    \underset{\nu\to+\infty}{\longrightarrow}
    \mathcal{M}_{u_\infty} (E)
    \mbox{ in the weak* topology of }
    \mathcal{D}' (I \times \reali^n;\reali)\,.
  \end{displaymath}
\end{theorem}

Above, we used the standard locally convex topologies and the related
notions. By boundedness in $\Lloc\infty (\reali^n;\reali)$ we mean
boundedness on every compact subset of $\reali^n$,
see~\cite[Definition~14.1]{MR2296978}. By convergence in
$\C1 (\reali;\reali)$ we mean the uniform convergence of the function
and of its derivative on every compact subset of $\reali$,
see~\cite[Example~I, Chapter~10]{MR2296978}.

A direct consequence of Theorem~\ref{thm:continuitynotx} worth
being noted is the following.

\begin{corollary}
  \label{cor:MayBeContNOTX}
  Let $f \in \mathscr{F}^n$,
  $g \in \C0 (\overline{I} \times \reali^n\times \reali;\reali)$,
  $u_o \in \L\infty (\reali; \reali)$ and
  $u \in \L\infty (I \times \reali^n; \reali)$.  $\mathcal{M}_u$ is
  continuous with respect to the locally convex topologies in
  $\C1 (\reali;\reali)$ and the weak* topology in
  $\mathcal{D}' (I \times \reali^n;\reali)$.
\end{corollary}

\begin{theorem}
  \label{thm:Identifynotx}
  Let $f \in \mathscr{F}^n$,
  $g \in \C0 (\overline{I} \times \reali^n\times \reali;\reali)$ and
  $u_o \in \L\infty (\reali; \reali)$. Assume that
  \begin{equation}
    \label{eq:33}
    \partial_u g \in
    \Lloc1\left(I; \L\infty (\reali^n\times\reali; \reali)\right) \,.
  \end{equation}
  Then,
  \begin{displaymath}
    \mbox{For all }
    u_1,u_2 \in \L\infty (I \times \reali^n; \reali)
    \quad
    \mathcal{M}_{u_1} = \mathcal{M}_{u_2}
    \quad \implies \quad
    u_1 = u_2 \,.
  \end{displaymath}
\end{theorem}

The following example shows the necessity of an assumption of
type~\eqref{eq:33} in Theorem~\ref{thm:Identifynotx}. Indeed, choose
$n=1$, $I = \reali_+$ and
\begin{displaymath}
  f (t,x,u) \colonequals x^2 \, u
  \,,\qquad  g (t,x,u) \colonequals 2 \, x \, u
  \qquad \mbox{ and } \qquad
  u_o (x) \colonequals 0\,,
\end{displaymath}
so that~\eqref{eq:33} fails and~\eqref{eq:1} becomes
\begin{equation}
  \label{eq:36}
  \left\{
    \begin{array}{l}
      \partial_t u + x^2 \, \partial_x u = 0
      \\
      u (0,x) = 0
    \end{array}
  \right.
\end{equation}
as long as smooth solutions are considered. A smooth solution
to~\eqref{eq:36} is $u_1 (t,x) \colonequals 0$. A further solution can
be defined as follows. Fix any function
\begin{displaymath}
  w \in \mathcal{D} (\reali; \reali)
  \quad \mbox{ with } \quad
  w (x) = 0 \iff x \leq 0 \,.
\end{displaymath}
Then, one can verify that also the map
\begin{displaymath}
  u_2 (t,x) \colonequals \left\{
    \begin{array}{l@{\qquad\mbox{if }}r@{\;}c@{\;}l}
      0
      & x
      & \geq
      & 0
      \\
      w\left(t + \frac1x\right)
      & x
      & <
      & 0
    \end{array}
  \right.
\end{displaymath}
is a strong solution to~\eqref{eq:36}. Since both $u_1$ and $u_2$ are
smooth solutions, by~\ref{item:33} we have
$\mathcal{M}_{u_1} = \mathcal{M}_{u_2}$, although $u_1 \neq u_2$.

With reference to~\eqref{eq:1}, the following quantity often recurs
below, where $(t,x,u) \in I \times \reali^n \times \reali$ and
$k \in \reali$:
\begin{equation}
  \label{eq:12}
  \Phi (t,x,u,k)
  \colonequals
  \sgn (u-k) \; \left(f (t,x,u) - f (t,x,k)\right)
\end{equation}
and we adopt the notation
\begin{equation}
  \label{eq:6}
  \sgn z \colonequals
  \left\{
    \begin{array}{r@{\qquad}r@{\,}c@{\,}l}
      1
      & z
      & >
      & 0
      \\
      0
      & z
      & =
      & 0
      \\
      - 1
      & z
      & <
      & 0 \,.
    \end{array}
  \right.
\end{equation}

\begin{definition}
  \label{def:muk}
  Let $f \in \mathscr{F}^n$,
  $g \in \C0 (\overline{I} \times \reali^n\times \reali;\reali)$,
  $u_o \in \L\infty (\reali; \reali)$ and
  $u \in \L\infty (I \times \reali^n; \reali)$. With reference
  to~\eqref{eq:1}, define for $k \in \reali$ the \emph{$k$-{Kru\v
      zkov} entropy production distribution}
  \begin{eqnarray*}
    \mu_k (\phi)
    & \colonequals
    & \int_{I} \int_{\reali^n}
      \left(
      \modulo{u (t,x) - k} \, \partial_t \phi (t,x)
      +
      \Phi\left(t, x, u (t,x),k\right) \cdot \grad_x \phi (t,x)
      \right) \d{x} \d{t}
    \\
    &
    &-
      \int_{I} \int_{\reali^n}
      \sgn\left(u (t,x) - k\right) \; \div_x f (t,x,k) \; \phi (t,x)
      \d{x} \d{t}
    \\
    &
    &+
      \int_{I} \int_{\reali^n}
      \sgn\left(u (t,x) - k\right) \; g \left(t,x,u (t,x)\right)
      \; \phi (t,x)
      \d{x} \d{t}
    \\
    &
    & + \int_{\reali^n} \modulo{u_o(x) - k} \, \phi(0,x) \d{x}
  \end{eqnarray*}
  where $\phi \in \mathcal{D} (I \times \reali^n;\reali)$ and we used
  the notation~\eqref{eq:12}, \eqref{eq:6}.
\end{definition}

The formal connection between $\mathcal{M}_u (E)$ and $\mu_k$ is
evident since, setting $E_k (u) = \modulo{u-k}$, we formally have
$\mu_k = \mathcal{M}_u (E_k)$. However, a rigorous extension of
Definition~\ref{def:muMnotx} to the case where $\partial_u E$ is
defined only up to sets of measure $0$ in $u$ is more delicate.  For
instance, with reference to~\ref{item:1notx} in
Proposition~\ref{prop:1notx}, consider two sequences $E_\nu$ and
$\check E_\nu$ such that
$E_\nu \underset{\nu\to+\infty}{\longrightarrow} E$,
$\check E_\nu \underset{\nu\to+\infty}{\longrightarrow} E$,
$E_\nu' \underset{\nu\to+\infty}{\longrightarrow} E'$,
$\check E_\nu' \underset{\nu\to+\infty}{\longrightarrow} E'$ pointwise
a.e. in $\reali$. Then, $\lim_{\nu\to+\infty} \mathcal{M}_u (E_\nu)$
and $\lim_{\nu\to+\infty} \mathcal{M}_u (\check E_\nu)$ might yield
different limits depending on the values of $E'$ on sets of measure
$0$. It is due to this point that we carefully specify the value of
$\sgn x$ in~\eqref{eq:6} for \emph{every} $x \in \reali$.

Nevertheless, an extension of Definition~\ref{def:muMnotx} to the case
$[u \mapsto E (u)] \in \AC (\reali; \reali)$, with $E'$ regulated, is
deferred to \S~\ref{subsec:case-class-entr}.

The properties of $\mu_k$ relevant to the sequel are proved in the
following proposition.

\begin{proposition}
  \label{prop:2}
  Let $f \in \mathscr{F}^n$,
  $g \in \C0 (\overline{I} \times \reali^n\times \reali;\reali)$,
  $u_o \in \L\infty (\reali; \reali)$ and
  $u \in \L\infty (I \times \reali^n; \reali)$. Then,
  \begin{enumerate}[label=\bf($\mathbf{\mu}$\arabic*)]
  \item \label{item:10} For all $k \in \reali$, $\mu_k$ is a
    distribution of order at most $1$.
  \item \label{item:11} For all compact sets $K_{tx}$ in
    $\reali \times \reali^n$ and all compact intervals $K_u$ in
    $\reali$, there exists a constant $C$ depending on $f,g,u$ --- but
    not on $u_o$ --- such that for all
    $\phi \in \mathcal{D} (I \times \reali^n; \reali)$ with
    $\spt\phi \subseteq K_{tx}$, the total variation of the map
    $k \mapsto \mu_k (\phi)$ on $K_u$ is bounded by
    \begin{displaymath}
      \tv \left(k \mapsto \mu_k (\phi); K_u\right)
      \leq C \, \norma{\phi}_{\C1 (I \times \reali^n; \reali)} \,.
    \end{displaymath}
  \item \label{item:27} For all $\overline{k} \in \reali$ and all
    $\phi \in \mathcal{D} (I \times \reali^n; \reali)$,
    \begin{displaymath}
      \mu_{\overline{k}} (\phi) = \dfrac12
      \left(
        \lim_{k \to \overline{k}-}\mu_k (\phi)
        +
        \lim_{k \to \overline{k}+}\mu_k (\phi)
      \right) \,.
    \end{displaymath}
  \end{enumerate}
\end{proposition}

\begin{remark}
  \label{rem:confusionMu}
  \rm In view of~\ref{item:10}, in the sequel we view $\mu_k$ also as
  a linear and continuous functional acting on
  $\Cc1 (I\times\reali^N;\reali)$, i.e.,
  $\mathcal{M}_u (E) \in \left(\Cc1 (I\times\reali^N;\reali)\right)'$.
\end{remark}

We now restrict our attention to maps $u$ that are distributional (or
weak) solutions to the Cauchy problem~\eqref{eq:1},
see~\cite[\S~4.3]{MR3468916}.

\begin{definition}
  \label{def:ws}
  Let $f \in \mathscr{F}^n$,
  $g \in \C0 (\overline{I} \times \reali^n\times \reali;\reali)$ and
  $u_o \in \L\infty (\reali; \reali)$. A function
  $u \in \L\infty (I \times \reali^n; \reali)$ is a
  \emph{distributional solution} to the Cauchy problem~\eqref{eq:1} if
  for all $\phi \in \Cc1 (I \times \reali^n; \reali)$
  \begin{eqnarray*}
    \int_{I} \int_{\reali^n}
    \left(
    u (t,x) \, \partial_t \phi (t,x)
    +
    f\left(t, x, u (t,x)\right) \cdot \grad_x \phi (t,x)
    \right) \d{x} \d{t}
    \\
    +
    \int_{I} \int_{\reali^n}
    g \left(t,x,u (t,x)\right) \; \phi (t,x)
    \d{x} \d{t}
    +\int_{\reali^n} u_o (x) \, \phi (0,x) \d{x}
    & =
    & 0 \,.
  \end{eqnarray*}
\end{definition}

\begin{remark}
  \label{rem:identity}
  \rm A direct consequence of Definition~\ref{def:muMnotx} and of
  Definition~\ref{def:ws} is that
  \begin{displaymath}
    u \mbox{ is a distributional solution to~\eqref{eq:1}}
    \iff
    \mathcal{M}_u (\Id_{\reali}) = 0 \,.
  \end{displaymath}
\end{remark}

\begin{remark}
  \label{rem:weakSol}
  \rm Knowing that $u$ is a distributional solution to~\eqref{eq:1}
  simplifies the above results thanks
  to~\cite[Lemma~1.3.3]{MR3468916}, which ensures that
  $\lim_{t\to0+}u (t) = u_o$ in the $\L\infty$ weak*
  topology. Therefore, in this case, $J_{u,u_o} = J_u$, with reference
  to~\eqref{eq:14}.
\end{remark}

\begin{remark}
  \label{rem:support}
  \rm Whenever $u$ is a distributional solution to~\eqref{eq:1}
  according to Definition~\ref{def:ws}, the map $k \mapsto \mu_k$ as
  introduced in Definition~\ref{def:muk} satisfies
  \begin{displaymath}
    k \leq \inf J_u
    \quad \mbox{ or } \quad
    k \geq \sup J_u
    \quad \implies \quad
    \mu_k = 0 \,.
  \end{displaymath}
\end{remark}

A further simplification is provided by the following proposition.

\begin{proposition}
  \label{prop:affinenotx}
  Let $f \in \mathscr{F}^n$,
  $g \in \C0 (\overline{I} \times \reali^n\times \reali;\reali)$ and
  $u_o \in \L\infty (\reali; \reali)$.
  \begin{enumerate}[\bf(\arabic*)]
  \item \label{item:30} Fix
    $u \in \L\infty (I \times \reali^n; \reali)$. Let
    $E \in \C1 (J_{u,u_o};\reali)$ and fix $b \in \reali$. Call
    \begin{displaymath}
      \tilde E (u) \colonequals E (u) + b \,.
    \end{displaymath}
    Then,
    \begin{displaymath}
      \mathcal{M}_u (E) = \mathcal{M}_u (\tilde E) \,.
    \end{displaymath}
  \item \label{item:29} Fix a distributional solution
    $u \in \L\infty (I \times \reali^n; \reali)$ to~\eqref{eq:1}. Let
    $E \in \C1 (J_u;\reali)$ and fix $a,b \in \reali$. Call
    \begin{displaymath}
      \tilde E (u) \colonequals E (u) + a \, u + b \,.
    \end{displaymath}
    Then,
    \begin{displaymath}
      \mathcal{M}_u (E) = \mathcal{M}_u (\tilde E) \,.
    \end{displaymath}
  \end{enumerate}
\end{proposition}

\noindent The above proposition shows that, in general,
$\mathcal{M}_u (E)$ may depend only on $E'$, as soon as $E$ is of
class $\C1$. On the contrary, when $u$ is a distributional solution,
$\mathcal{M}_u (E)$ may depend only on $E''$, as soon as $E$ is of
class $\C2$. Indeed, the next theorem provides an explicit formula
that, in particular, also confirms the last part of this observation.

\begin{theorem}
  \label{thm:formulaC2}
  Let $f \in \mathscr{F}^n$,
  $g \in \C0 (\overline{I} \times \reali^n\times \reali;\reali)$,
  $u_o \in \L\infty (\reali;\reali)$ and
  $u \in \L\infty (I \times \reali^n; \reali)$. Then, the following
  statements are equivalent:
  \begin{enumerate}[\bf(\arabic*)]
  \item \label{item:21} $u$ is a distributional solution
    to~\eqref{eq:1} in the sense of Definition~\ref{def:ws};
  \item \label{item:22} For all $E \in \C2 (\reali;\reali)$, with the
    notation in Definition~\ref{def:muMnotx} and in
    Definition~\ref{def:muk}, the following formula holds:
    \begin{equation}
      \label{eq:2}
      \mathcal{M}_u (E)
      =
      \dfrac12 \,
      \int_{J_u} E'' (k) \; \mu_k \d{k} \,.
    \end{equation}
  \end{enumerate}
\end{theorem}

\noindent Formula~\eqref{eq:2} is an equality between
distributions. In fact, in view of~\ref{item:4notx} in
Proposition~\ref{prop:1notx}, by~\eqref{eq:2} the action of
$\mathcal{M}_u (E)$ on a test function $\phi$ is
\begin{displaymath}
  \left(\mathcal{M}_u (E)\right) (\phi)
  =
  \dfrac12 \,
  \int_{J_u} E'' (k) \;
  \mu_k \left(\phi\right)    \d{k}\,.
\end{displaymath}
On the other hand, Proposition~\ref{prop:integral} can be applied with
\begin{displaymath}
  d=n+1\,,\quad
  \Omega = I \times \reali^n \,,\quad
  (X,\mathcal{S},\mathfrak{m}) = (J_u, \mathcal{B}, \mathcal{L})\,,\quad
  \Upsilon \colon  k \mapsto E'' (k) \, \mu_k\,,
\end{displaymath}
ensuring that the integral $\int_{\reali} E'' (k) \; \mu_k \d{k}$ is
well defined and a distribution of order $1$ on $I \times \reali^n$.

\subsection{Extension to \texorpdfstring{$E'$}{E'} Regulated}
\label{subsec:case-class-entr}

Our next goal is to provide a unique definition comprising both
Definition~\ref{def:muMnotx} and Definition~\ref{def:muk}, while
keeping the validity of the representation formula~\eqref{eq:2}. To
this aim, we adopt the standard notation
$\psi (x-) \colonequals \lim_{\xi \to x-} \psi (\xi)$,
$\psi (x+) \colonequals \lim_{\xi \to x+} \psi (\xi)$ and recall the
following definition.

\begin{definition}[{\cite[Definition~4.1.1]{MR3839599}}]
  \label{def:regulated}
  Fix a map $\Psi \colon \reali \to \reali$. $\Psi$ is
  \emph{regulated} if for all $x \in \reali$, both limits $\Psi (x-)$
  and $\Psi (x+)$ exist and are finite.
\end{definition}

By~\cite[Theorem~4.1.8]{MR3839599}, if $\Psi$ is regulated, by
modifying $\Psi$ on at most countable points, we may assume that for
all $x \in \reali$,
$\Psi (x) = \frac12 \left(\Psi (x-) + \Psi (x+)\right)$. Thus, we now
reduce the regularity assumption on $E$ requiring that $E \in \ACR (\reali;\reali)$, where
\begin{equation}
  \label{eq:20}
  \ACR (\reali; \reali)
  \colonequals
  \left\{
    E \in \AC (\reali; \reali)
    \colon
    \begin{array}{@{}c@{}}
      E' \mbox{ admits a regulated representative satisfying }
      \\
      \forall \,u \in \reali \qquad E' (u) = \dfrac12 \left( E' (u-) + E' (u+) \right)
    \end{array}
  \right\}
\end{equation}

\noindent The choice~\eqref{eq:20} is coherent with~\eqref{eq:6} and
this allows to extend Definition~\ref{def:muMnotx} introducing
\begin{equation}
  \label{eq:37}
  \overline{\mathcal{M}}_u
  \colon
  \ACR (\reali; \reali)
  \to
  \mathcal{D}' (I \times \reali^n;\reali)
\end{equation}
by means of the same formal expression~\eqref{eq:21}.  Note that
$\overline{\mathcal{M}}_u$ extends Definition~\ref{def:muMnotx} of
$\mathcal{M}_u$ so that it comprises the Definition~\ref{def:muk} of
$\mu_k$. Below, for simplicity, we write $\mathcal{M}_u$ for
$\overline{\mathcal{M}}_u$ and it results --- now not only formally
--- that
\begin{displaymath}
  \mu_k = \overline{\mathcal{M}}_u (E_k)
  \quad \mbox{ where } \quad
  E_k (u) \colonequals \modulo{u-k} \,.
\end{displaymath}

\begin{remark}
  \label{rem:ACR}
  \rm $\ACR (\reali;\reali)$, with the topology induced by the uniform
  convergence of $E$ and $E'$ on compact subsets of $\reali$, is a
  Fr\'echet space, see~\cite[Chapter~4]{MR3839599}.
\end{remark}

The following result is thus of interest and of use in the sequel.

\begin{proposition}
  \label{prop:regCont}
  Let $f \in \mathscr{F}^n$,
  $g \in \C0 (\overline{I} \times \reali^n\times \reali;\reali)$,
  $u_o \in \L\infty (\reali; \reali)$ and
  $u \in \L\infty (I \times \reali^n; \reali)$. For all
  $\nu \in \naturali$, let $E,E_\nu \in \ACR (\reali;\reali)$ as
  in~\eqref{eq:20}. If
  \begin{displaymath}
    E_\nu \underset{\nu\to+\infty}{\longrightarrow} E
    \mbox{ and }
    E'_\nu \underset{\nu\to+\infty}{\longrightarrow} E'
    \mbox{ uniformly on } J_{u,u_o}
  \end{displaymath}
  then for any $\phi \in \Cc1 (I\times\reali^n;\reali)$,
  $\left( \overline{\mathcal{M}_u} (E_\nu)\right) (\phi)
  \underset{\nu\to+\infty}{\longrightarrow} \left(
    \overline{\mathcal{M}_u} (E)\right) (\phi)$.
\end{proposition}

The next result extends the representation formula~\eqref{eq:2} from
the case of Definition~\ref{def:muMnotx} to that of~\eqref{eq:37} on
the basis of the formal identity $\d E' (k) = E'' (k) \, \d{k}$.

\begin{theorem}
  \label{thm:formulaC1}
  Let $f \in \mathscr{F}^n$,
  $g \in \C0 (I \times \reali^n \times \reali; \reali)$ and
  $u_o \in \L\infty (\reali; \reali)$. Fix a distributional solution
  $u \in \L\infty (\reali;\reali)$ to~\eqref{eq:1} in the sense of
  Definition~\ref{def:ws}. If $E \in \ACR (\reali;\reali)$ as
  in~\eqref{eq:20}, then with the notation in~\eqref{eq:21} and in
  Definition~\ref{def:muk},
  \begin{equation}
    \label{eq:3}
    \overline{\mathcal{M}}_u (E)
    =
    \dfrac12 \, \int_{\reali} \mu_k \, \d{E' (k)} \,.
  \end{equation}
\end{theorem}

Recall that $k \mapsto \mu_k (\phi)$ is in $\BV (\reali;\reali)$
by~\ref{item:11} in Proposition~\ref{prop:2}, supported in $J_u$ by
Remark~\ref{rem:support}. Thus, Proposition~\ref{prop:integral2}
ensures that the Kurzweil-Stieltjes integral
$\int_{\reali} \mu_k \d{E' (k)}$ in~\eqref{eq:3} is well defined and a
distribution. Moreover, this integral requires $E'$ to be regulated,
but its value is independent from the choice~\eqref{eq:20}. On the
contrary, the left hand side in~\eqref{eq:3} does depend on the
specific values of $E'$ at its points of jump. This apparent
contradiction vanishes when considering the role of our particular
choice~\eqref{eq:6} that specifies $\mu_k$ in
Definition~\ref{def:muk}.

More precisely, it is reasonable to expect that a coherent replacement
of~\eqref{eq:6} and~\eqref{eq:20} with
\begin{displaymath}
  \sgn z \colonequals
  \left\{
    \begin{array}{l@{\qquad}r@{\,}c@{\,}l}
      1
      & z
      & >
      & 0
      \\
      1-2\,\theta
      & z
      & =
      & 0
      \\
      - 1
      & z
      & <
      & 0 \,.
    \end{array}
  \right.
  \quad \mbox{ and } \qquad
  E' (u) = \theta \; E' (u-) + (1-\theta) \, E' (u+)
\end{displaymath}
for $\theta \in [0,1]$ might keep the validity of~\eqref{eq:3}. This
proof relies on an extension of~\cite[(iii) in
Exercise~4.2.5]{MR3839599} where regularized functions are defined
through left limits, weighted averages or right limits.

Remark that under the stronger regularity
$E \in \C2 (\reali; \reali)$, the equality~\eqref{eq:2} directly
follows from~\eqref{eq:3} by~\cite[Corollary~5.4.4 and
Theorem~6.2.12]{MR3839599}.

\begin{remark}
  \rm Note that a naive approach to the proof of~\eqref{eq:3} might
  rely on~\eqref{eq:2}, then approximating a regulated $E'$ with a
  sequence of $E'_n$ in $\C1$. However, $\C1$ is not dense in the
  space of regulated functions and limits in~\eqref{eq:3} on $E'$
  require uniform convergence, see~\cite[(i) in
  Theorem~6.3.9]{MR3839599}.
\end{remark}

\begin{corollary}
  \label{cor:solution}
  Let $f \in \mathscr{F}^n$,
  $g \in \C0 (\overline{I} \times \reali^n \times \reali; \reali)$ and
  $u_o \in \L\infty (\reali; \reali)$. Fix a distributional solution
  $u \in \L\infty (I\times\reali^N; \reali)$ to~\eqref{eq:1} in the
  sense of Definition~\ref{def:ws}. Then, the following statements are
  equivalent:
  \begin{enumerate}[\bf(\arabic*)]
  \item \label{item:24} For all $k \in \reali$, $\mu_k$ is a non
    negative Radon measure.
  \item \label{item:25} For every convex entropy $E$ of class
    $\C\infty$, $\mathcal{M}_{u} (E)$ is a non negative Radon measure.
  \item \label{item:26} Every convex entropy $E$ is in
    $\ACR (\reali;\reali)$ as in~\eqref{eq:20} and
    $\overline{\mathcal{M}}_{u} (E)$ is a non negative Radon measure.
  \end{enumerate}
\end{corollary}

\noindent If $n=1$, $g\equiv 0$ and $f = f (x,u)$,
Condition~\ref{item:24} in Corollary~\ref{cor:solution} guarantees the
uniqueness of $u$ as entropy solution to~\eqref{eq:1},
see~\cite[Theorem~2.6]{ColomboPerrollazDepX}. In the general case,
\cite[Theorem~1]{Kruzkov1970} requires~\ref{item:24} in
Corollary~\ref{cor:solution} together with a trace condition for
$t \to 0+$ to ensure uniqueness. To bypass this requirement, we refer
to the discussions in~\cite[\S~4.5]{MR3468916}, as well as
in~\cite{MR1771520}, \cite{MR1235201}, \cite[Chapter~4]{MR2064166}
and~\cite{MR1869441}.

The introduction of $\mathcal{M}_u$ in Definition~\ref{def:muMnotx}
--- and of its extended version $\overline{\mathcal{M}}_u$
in~\eqref{eq:37} -- allows a new formulation of known results
concerning the uniqueness of solutions to scalar one dimensional
conservation laws. Examples are the following corollaries: for their
proofs we refer to the original ones as detailed below.

\begin{corollary}[{\cite[Theorem~1]{MR1296003},
    \cite[Corollary~2.5]{MR2104269}}]
  \label{cor:Panov}
  Let $f \in \C1 (\reali; \reali)$ be convex and
  $u_o \in \L\infty (\reali^n;\reali)$. Fix a distributional solution
  $u \in \L\infty (I \times\reali^N; \reali)$ to
  \begin{equation}
    \label{eq:22}
    \left\{
      \begin{array}{l}
        \partial_t u + \partial_x f (u) =0
        \\
        u (0,x) = u_o (x)
      \end{array}
    \right.
  \end{equation}
  in the sense of Definition~\ref{def:ws}. If there exists
  $E_* \in \C2 (\reali; \reali)$ uniformly convex and such that
  $\mathcal{M}_u (E_*)$ is non negative, then
  \begin{enumerate}[label=\bf(\arabic*)]
  \item \label{item:2} For all $k \in \reali$, $\mu_k$ is non
    negative.
  \item \label{item:3} For all $E \in \C2 (\reali; \reali)$ convex,
    $\mathcal{M}_u (E)$ is non negative.
  \end{enumerate}
\end{corollary}

\begin{corollary}[{\cite[Theorem~1]{Kruzkov1970},
    \cite[Definition~3]{MR2592291}}]
  \label{cor:Kruzkov}
  Let $f \in \mathscr{F}^n$,
  $g \in \C0 (\overline{I} \times \reali^n \times \reali; \reali)$ and
  $u_o \in \L\infty (\reali; \reali)$. Assume that
  \begin{equation}
    \label{eq:4}
    \forall\, (t,x) \in \reali_+\times\reali^n
    \quad
    \forall\, \mathbf{n} \in \reali^n \setminus\{0\}
    \qquad
    f (t,x,u) \cdot \mathbf{n} \neq 0
    \quad \mbox{ for a.e. } u \in \reali \,.
  \end{equation}
  If $u_1,u_2 \in \L\infty (I \times \reali^n; \reali)$ are
  distributional solutions to~\eqref{eq:1} such that for all convex
  $E \in \C2 (\reali; \reali)$ it holds that
  $\mathcal{M}_{u_1} (E) \geq 0$ and $\mathcal{M}_{u_2} (E) \geq 0$,
  then $u_1 = u_2$.
\end{corollary}

The proof follows from the classical {Kru\v zkov} uniqueness result
in~\cite[Theorem~1]{Kruzkov1970} and from the discussion
in~\cite[Definition~3]{MR2592291} which ensures that the trace at $0$
condition in~\cite[Definition~1]{Kruzkov1970} is not necessary under
assumption~\eqref{eq:4}.

Remark that in the particular case $n=1$, $f = f (x,u)$ and
$g\equiv 0$, \cite[Theorem~2.6]{ColomboPerrollazDepX} shows that the
genuine non linearity condition~\eqref{eq:4} is not necessary.

When the two corollaries above can be applied together, they yield the
uniqueness of solutions to~\eqref{eq:1}. In the case of systems, this
procedure acquires relevance due to the well known lack of
entropies. Refer to~\cite{MR4487515} for a result in this direction
and to~\cite{MR4661213} for a sort of opposite point of view, not
explicitly relying on entropies.

\subsection{Fourier Transform}
\label{subsec:fourier-transform}

All through this section, we consider complex valued entropies, i.e.,
$E \in \mathscr{F} + i \mathscr{F}$, so that the corresponding entropy
flux is $F \in \mathscr{F}^n + i \mathscr{F}^n$. At the same time,
also test functions are assumed to be complex valued
$\phi \in \mathcal{D} (I\times\reali^N; \complessi)$. The extension of
the results in the previous sections to this setting is immediate,
thanks to the linearity of expression defining $\mathcal{M}_u$ in both
$E$ and $\phi$.

Fix a distributional solution $u$ to~\eqref{eq:1} in the sense of
Definition~\ref{def:ws}. Definition~\ref{def:muk} introduces the map
$\mu$, which we now see as a function
\begin{displaymath}
  \begin{array}{ccccc}
    \mu
    & \colon
    & \reali
    & \to
    & \mathcal{D}' (I \times \reali^n; \complessi)
    \\
    &
    & k
    & \mapsto
    & \mu_k \,.
  \end{array}
\end{displaymath}
By Remark~\ref{rem:support}, the map $k \mapsto \mu_k$ is compactly
supported and its Fourier transform~\cite[\S~VII.1]{MR1065993}
\begin{equation}
  \label{eq:15}
  \widehat\mu (\xi)
  \colonequals
  \int_{\reali} \mu_k \; e^{-i\, \xi \, k} \d{k}
\end{equation}
is well defined for all $\xi \in \complessi$ thanks to
Proposition~\ref{prop:integral} and provides relevant information.

\begin{proposition}
  \label{prop:Fourier}
  Let $f \in \mathscr{F}^n$,
  $g \in \C0 (\overline{I} \times \reali^n \times \reali; \reali)$ and
  $u_o \in \L\infty (\reali; \reali)$. Let
  $u \in \L\infty (I\times\reali^N; \reali)$ be a distributional
  solution to~\eqref{eq:1} in the sense of Definition~\ref{def:ws}.
  Then,
  \begin{enumerate}[label=\bf(FT.\arabic*)]
  \item \label{item:19} $\widehat\mu$ as defined in~\eqref{eq:15} is
    an entire function defined on $\complessi$ with values in
    $\mathcal{D}' (I \times \reali^n; \complessi)$.
  \item \label{item:20} Setting for $\xi \in \complessi$ and
    $k \in \reali$
    \begin{displaymath}
      \mathcal{E}_{\xi} (k)
      \colonequals
      - 2 \, \dfrac{e^{-i\, \xi \, k}-1+i\, \xi \, k}{\xi^2} \,,
    \end{displaymath}
    with the notation in Definition~\ref{def:muMnotx} and
    in~\eqref{eq:15}, the following equality holds:
    \begin{displaymath}
      \widehat\mu (\xi)
      =
      \mathcal{M}_u (\mathcal{E}_\xi) \,.
    \end{displaymath}
  \end{enumerate}
\end{proposition}

\noindent Inspired by the {Kru\v zkov} entropy
$E_k (w) = \modulo{w-k}$, we consider a one parameter family of
entropies
\begin{displaymath}
  \begin{array}{ccccc}
    E
    & \colon
    & \reali \times \reali
    & \to
    & \reali
    \\
    &
    & (w,k)
    & \mapsto
    & E (w,k)
  \end{array}
  \quad \mbox{ with flux } \quad
  \partial_w F (t,x,w,k) = \partial_w E (w,k) \; \partial_w f (t,x,w) \,.
\end{displaymath}
Here, $k$ plays the role of a parameter with the same dimension as $u$
and the roles of $k$ and $w$ can occasionally be exchanged.

Given $u \in \L\infty (I \times \reali: \reali)$, interpret
Definition~\ref{def:muMnotx} of $\mathcal{M}_u$ as follows:
\begin{equation}
  \label{eq:18}
  \begin{array}{ccccl}
    \mathcal{M}_u (E)
    & \colon
    & \reali
    & \mapsto
    & \mathcal{D}' (I \times \reali^n; \complessi)
    \\
    &
    & k
    & \mapsto
    & \mathcal{M}_u \left(E (\cdot,k)\right) \,.
  \end{array}
\end{equation}
and recall that if $u$ is a distributional solution to~\eqref{eq:1},
then by Theorem~\ref{thm:formulaC2}
\begin{equation}
  \label{eq:38}
  \mathcal{M}_u \left(E (\cdot,k)\right)
  =
  \dfrac12 \int_{J_{u,u_o}} \partial^2_{ww} E (w,k) \, \mu_w \d{w} \,.
\end{equation}
Then, by the linearity and continuity proved in~\ref{item:1notx} in
Proposition~\ref{prop:1notx},
\begin{displaymath}
  \left(\widehat{\mathcal{M}_u (E)}\right) (\xi)
  =
  \mathcal{M}_u \left(\widehat{E} (\cdot ,\xi)\right)
\end{displaymath}
where
\begin{displaymath}
  \left(\widehat{{\mathcal{M}_u (E)}}\right) (\xi)
  =
  \int_{\reali} e^{-i\,k\,\xi} \;
  \mathcal{M}_u \left(E (\cdot,k)\right) \,
  \d{k}
  \quad \mbox{ and } \quad
  \widehat{E} (\cdot ,\xi)
  =
  \int_{\reali} e^{-i\,k\,\xi} \; E (\cdot,k) \, \d{k} \,.
\end{displaymath}

Define below the translation operator
\begin{equation}
  \label{eq:23}
  \left(\mathcal{T}_k E\right) (u) = E (u-k) \,.
\end{equation}
Observe that, with reference to {Kru\v zkov} entropies
$E_k (u) = \modulo{u-k}$, clearly $E_k = \mathcal{T}_kE_0$.

\begin{theorem}
  \label{thm:propFT}
  Let $f \in \mathscr{F}^n$ and
  $g \in \C0 (\overline{I} \times \reali^n \times \reali; \reali)$.
  Consider an entropy $E \in \C2 (\reali;\reali)$ such that $E''$ is a
  tempered distribution and $\spt \widehat {E''}$ admits a cluster
  point.  For any $u_o,v_o \in \L\infty (\reali^n;\reali)$ consider
  the Cauchy problems
  \begin{equation}
    \label{eq:39}
    \left\{
      \begin{array}{l}
        \partial_t u + \div_x f (t,x,u) = g (t,x,u)
        \\
        u (0,x) = u_o (x) \,.
      \end{array}
    \right.
    \quad \mbox{ and } \quad
    \left\{
      \begin{array}{l}
        \partial_t v + \div_x f (t,x,v) = g (t,x,v)
        \\
        v (0,x) = v_o (x) \,.
      \end{array}
    \right.
  \end{equation}
  For any distributional solutions
  $u,v \in \L\infty (I \times \reali^n;\reali)$ to~\eqref{eq:39}, call
  $\mathcal{M}^{u_o}_u$, $\mathcal{M}^{v_o}_v$ the corresponding
  entropy production functionals as in
  Definition~\ref{def:muMnotx}. Then, with the notation~\eqref{eq:23}.
  \begin{displaymath}
    \forall \,k \in \reali \qquad
    \mathcal{M}^{u_o}_u \left(\mathcal{T}_kE \right) =
    \mathcal{M}^{v_o}_v \left(\mathcal{T}_kE\right) \quad \implies \quad
    \mathcal{M}^{u_o}_u = \mathcal{M}^{v_o}_v \,,
  \end{displaymath}
\end{theorem}

\noindent For the notion of tempered distribution, refer
to~\cite[Definition~7.1.7]{MR1065993}.

If also~\eqref{eq:33} holds and $u_o = v_o$, then
Theorem~\ref{thm:propFT} also ensures that $u = v$ thanks to Theorem~\ref{thm:Identify}.

\section{Time and Space Dependent Entropy Operator}
\label{sec:time-space-dependent}

With reference to~\eqref{eq:1}, we assume that the flux $f$ is in
$\mathscr{F}^n$, the source $g$ is in
$\C0 (\overline{I} \times \reali^n \times \reali; \reali)$ and the
initial datum $u_o$ is in $\L\infty (\reali; \reali)$.

Here we extend Definition~\ref{def:eefp} to the case of $t$ and $x$
dependent \emph{``entropies''}.

\begin{definition}
  \label{def:EEF}
  Fix a non trivial real interval $J$.  A pair of functions $(E,F)$ is
  an \emph{entropy -- entropy flux pair} with respect to $f$ if for
  all $t \in I$, $x \in \reali^n$
  \begin{eqnarray}
    \nonumber
    u \mapsto E (t,x,u)
    & \in
    & \AC (\reali; \reali)\,;
    \\
    \nonumber
    u \mapsto F (t,x,u)
    & \in
    & \AC (\reali; \reali^n) \,;
    \\
    \label{eq:75}
    \mbox{for a.e. } u \in J \qquad
    \partial_u F (t, x,u)
    & =
    & \partial_u E (t,x,u) \;\; \partial_u f (t, x,u) \,.
  \end{eqnarray}
\end{definition}

By the chain rule, as soon as a strong, or classical, solution $u$
to~\eqref{eq:1} is available, an entropy-entropy flux pair $(E,F)$ of
class $\C1$ satisfies for all $(t,x) \in I \times \reali^n$
\begin{eqnarray}
  \nonumber
  &
  & \partial_t \left(E\left(t,x,u (t,x)\right)\right)
    +
    \div_x \left(F\left(t,x,u (t,x)\right)\right)
  \\
  \label{eq:13}
  & =
  &(\partial_t E) \left(t,x,u (t,x)\right)
    +
    (\div_x F) \left(t,x,u (t,x)\right)
  \\
  \nonumber
  &
  & + \left(g\left(t,x,u (t,x)\right)
    -
    (\div_x f)\left(t,x,u (t,x)\right)\right) \,
    (\partial_u E) \left(t,x,u (t,x)\right) \,.
\end{eqnarray}

We proceed extending also Definition~\ref{def:muMnotx}.

\begin{definition}
  \label{def:muM}
  Let $f \in \mathscr{F}^n$,
  $g \in \C0 (\overline{I}\times \reali^n\times \reali;\reali)$,
  $u_o \in \L\infty (\reali; \reali)$ and
  $u \in \L\infty (I \times \reali^n; \reali)$. Define the
  \emph{entropy production distribution}
  \begin{displaymath}
    \begin{array}{ccccc}
      \mathcal{M}_u^{tx}
      & \colon
      & \mathscr{F}_c
      & \to
      & \reali
      \\
      &
      & E
      & \mapsto
      & \mathcal{M}_u^{tx} (E)
    \end{array}
  \end{displaymath}
  where for $F \in \mathscr{F}_c^n$ satisfying~\eqref{eq:75},
  \begin{eqnarray}
    \nonumber
    \mathcal{M}_u^{tx} (E)
    & \colonequals
    & -
      \int_{I} \int_{\reali^n}
      (\partial_u E)\left(t, x, u (t,x)\right)
      \; (\div_x f)\left(t, x, u (t,x)\right)
      \d{x} \d{t}
    \\
    \label{eq:40}
    &
    & +
      \int_{I} \int_{\reali^n}
      \left(
      (\div_x F)\left(t, x, u (t,x)\right)
      +
      (\partial_t E)\left(t, x, u (t,x)\right)
      \right)
      \d{x} \d{t}
    \\
    \nonumber
    &
    & +
      \int_{I} \int_{\reali^n}
      (\partial_u E)\left(t, x, u (t,x)\right) \;
      g \left(t,x,u (t,x)\right)
      \d{x} \d{t}
    \\
    \nonumber
    &
    & +\int_{\reali^n} E\left(0, x, u_o (x)\right) \d{x} \,.
  \end{eqnarray}
\end{definition}
In the above definition, $E$ is not required to be compactly supported
in $u$ since $u$ is assumed throughout to be bounded. An alternative
construction may consist in restricting $E$ to be compactly supported
in all variables, while allowing $u$ to be merely measurable.

\begin{remark}
  \label{rem:identitytx}
  \rm A direct consequence of Definition~\ref{def:ws} and of
  Definition~\ref{def:muM} is that
  \begin{displaymath}
    u \mbox{ is a distributional solution to~\eqref{eq:1}}
    \iff
    \mathcal{M}^{tx}_u (\phi \otimes \Id_{\reali}) = 0
    \quad \mbox{ for all } \phi \in \Cc1 (I\times\reali^n;\reali) \,.
  \end{displaymath}
\end{remark}

\begin{lemma}
  \label{lem:extension}
  Let $f \in \mathscr{F}^n$,
  $g \in \C0 (\overline{I}\times \reali^n\times \reali;\reali)$,
  $u_o \in \L\infty (\reali; \reali)$ and
  $u \in \L\infty (I \times \reali^n; \reali)$. For all
  $E \in \C1 (\reali;\reali)$ and all
  $\phi \in \Cc1 (I\times\reali^n;\reali)$,
  \begin{equation}
    \label{eq:47}
    \left(\mathcal{M}_u (E)\right) (\phi)
    =
    \mathcal{M}^{tx}_u (\phi \otimes E)  \,.
  \end{equation}
\end{lemma}

\noindent The proof follows by direct computations, using
Definition~\ref{def:muMnotx} and Definition~\ref{def:muM}.

We now provide the basic properties of $\mathcal{M}^{tx}$ and, in
particular, we verify that $\mathcal{M}_u^{tx} (E)$ is well defined
only as a function of $E$ and is independent of the particular choice
of the flux $F$, as long as $(E,F)$ is an entropy - entropy flux pair
in the sense of Definition~\ref{def:EEF}.

\begin{proposition}
  \label{prop:1}
  Let $f \in \mathscr{F}^n$,
  $g \in \C0 (\overline{I}\times \reali^n\times \reali;\reali)$,
  $u_o \in \L\infty (\reali; \reali)$ and
  $u \in \L\infty (I \times \reali^n; \reali)$. Then,
  \begin{enumerate}[label=\bf($\mathcal{\mathcal{M}}^{tx}$\arabic*)]
  \item \label{item:9} Let $F$ and $\check F$ be two entropy fluxes
    related to $E$ as in~\eqref{eq:75} and call ${\mathcal{M}}_u^{tx}$
    and $\check{\mathcal{M}}_u^{tx}$ the corresponding functions
    resulting from Definition~\ref{def:muM}. Then,
    ${\mathcal{M}}_u^{tx} = \check{\mathcal{M}}_u^{tx}$.
  \item \label{item:34} If $u$ is a strong solution to~\eqref{eq:1}
    with initial datum $u_o$, then $\mathcal{M}_u^{tx} \equiv 0$.
  \item \label{item:1} $\mathcal{M}_u^{tx}$ is linear and if
    $\rest{E}{I \times \reali^n \times J_{u,u_o}} = 0$, then
    $\mathcal{M}_u^{tx} (E) =0$.
  \end{enumerate}
\end{proposition}

\noindent In the spirit of Remark~\ref{rem:Ispt}, an alternative
formulation of the latter part of~\ref{item:1} is, for arbitrary
$E_1, E_2 \in \mathscr{F}_c$,
\begin{displaymath}
  \mbox{For all } (t,x,u) \in I \times \reali^n \times J_{u,u_o}
  \quad E_1 (t,x,u) = E_2 (t,x,u)
  \quad \implies \quad
  \mathcal{M}_u (E_1) = \mathcal{M}_u (E_2) \,.
\end{displaymath}
Alternatively, one can say that $\mathcal{M}_u$ is \emph{supported} in
$\C1 (I \times \reali^n \times J_{u,u_o}; \reali) \cap \mathscr{F}_c$.

\begin{theorem}
  \label{thm:continuity}
  For $\nu \in \naturali$, let $f,f_\nu \in \mathscr{F}^n$;
  $g, g_\nu \in \C0 (\overline{I} \times \reali^n \times
  \reali;\reali)$; $u_o, u_o^\nu \in \L\infty (\reali; \reali)$ and
  $u_\infty, u^\nu \in \L\infty (I \times \reali^n; \reali)$;
  $E,E_\nu \in \mathscr{F}_c$. Assume that
  \begin{eqnarray*}
    f_\nu
    & \underset{\nu\to+\infty}{\longrightarrow}
    & f \mbox{ in } \mathscr{F}^n
      \mbox{ with respect to~\eqref{eq:45}};
    \\
    g_\nu
    & \underset{\nu\to+\infty}{\longrightarrow}
    & g \mbox{ uniformly on compact subsets of }
      I \times \reali^n \times \reali;
    \\
    u_o^\nu
    & \underset{\nu\to+\infty}{\longrightarrow}
    & u_o \mbox{ pointwise a.e.~in } \reali^n \mbox{ and }
      u_o^\nu \mbox{ is bounded in } \Lloc\infty (\reali^n; \reali);
    \\
    u^\nu
    & \underset{\nu\to+\infty}{\longrightarrow}
    & u_\infty \mbox{ pointwise a.e.~in } \reali \times \reali^n \mbox{ and }
      u^\nu \mbox{ is bounded in } \Lloc\infty (I \times \reali^n; \reali);
    \\
    E_\nu
    & \underset{\nu\to+\infty}{\longrightarrow}
    & E \mbox{ in } \mbox{ in } \mathscr{F}_c
      \mbox{ with respect to the inductive topology}.
  \end{eqnarray*}
  \noindent Call $\mathcal{M}^{tx}_{u_\infty}$ the entropy production
  functional introduced in Definition~\ref{def:muM} and
  $\mathcal{M}^{tx}_\nu$ the entropy production functional generated
  by the Cauchy problem
  \begin{equation}
    \label{eq:48}
    \left\{
      \begin{array}{l}
        \partial_t u + \div_x f_\nu (t,x,u) = g_\nu (t,x,u)
        \\
        u (0,x) = u_o^\nu (x)
      \end{array}
    \right.
  \end{equation}
  and by the function $u_\nu$, according to
  Definition~\ref{def:muM}. Then,
  $\mathcal{M}^{tx}_\nu (E_\nu)
  \underset{\nu\to+\infty}{\longrightarrow}
  \mathcal{M}^{tx}_{u_\infty} (E)$.
\end{theorem}

We stress that in the above statement, neither $u_\infty$ nor any of
the $u_\nu$ are required to be in any sense solution to the Cauchy
problems~\eqref{eq:1} and~\eqref{eq:48}.

For a characterization of the convergence
$E_\nu \underset{\nu\to+\infty}{\longrightarrow} E$ with respect to
the inductive topology in $\mathscr{F}_c$, refer to
Lemma~\ref{lem:treves}.

Note that our choice of having weak requirements on the convergences
$u_o^\nu \underset{\nu\to+\infty}{\longrightarrow} u_o$ and
$u^\nu \underset{\nu\to+\infty}{\longrightarrow} u_\infty$ prevents to
rephrase Theorem~\ref{thm:continuity} as a continuity result
referred to suitable topologies.

Again, we underline the following direct consequence of
Theorem~\ref{thm:continuity}.

\begin{corollary}
  \label{cor:MayBeCont}
  $\mathcal{M}_u^{tx}$ is continuous with respect to the inductive
  topology in $\mathscr{F}_c$.
\end{corollary}

\begin{theorem}
  \label{thm:Identify}
  Let $f \in \mathscr{F}^n$,
  $g \in \C0 (\overline{I}\times \reali^n\times \reali;\reali)$ and
  $u_o \in \L\infty (\reali; \reali)$. Assume that~\eqref{eq:33}
  holds. Then,
  \begin{displaymath}
    \mbox{For all }
    u_1,u_2 \in \L\infty (I \times \reali^n; \reali) \quad
    \mathcal{M}_{u_1}^{tx} = \mathcal{M}_{u_2}^{tx}
    \quad \implies \quad
    u_1 = u_2 \,.
  \end{displaymath}
\end{theorem}

\begin{proposition}
  \label{prop:affine}
  Let $f \in \mathscr{F}^n$,
  $g \in \C0 (\overline{I}\times \reali^n\times \reali;\reali)$ and
  $u_o \in \L\infty (\reali^n; \reali)$.
  \begin{enumerate}[\bf(\arabic*)]
  \item \label{item:31} Fix
    $u \in \L\infty (I \times \reali^n; \reali)$. Let
    $E \in \mathscr{F}$ and fix
    $b \in \Cc1 (I \times \reali^n; \reali)$. Call
    \begin{displaymath}
      \tilde E (t,x,u) = E (t,x,u) + b (t,x) \,.
    \end{displaymath}
    Then,
    \begin{displaymath}
      \mathcal{M}_u^{tx} (E) = \mathcal{M}_u^{tx} (\tilde E) \,.
    \end{displaymath}
  \item \label{item:32} Fix a distributional solution
    $u \in \L\infty (I \times \reali^n; \reali)$. Let
    $E \in \mathscr{F}$ and fix
    $a,b \in \Cc1 (I \times \reali^n; \reali)$. Call
    \begin{displaymath}
      \tilde E (u) = E (t,x,u) + a (t,x) \, u + b (t,x) \,.
    \end{displaymath}
    Then,
    \begin{displaymath}
      \mathcal{M}_u^{tx} (E) = \mathcal{M}_u^{tx} (\tilde E) \,.
    \end{displaymath}
  \end{enumerate}
\end{proposition}

\begin{theorem}
  \label{thm:formulaC2exProp}
  Let $f \in \mathscr{F}^n$,
  $g \in \C0 (\overline{I} \times \reali^n\times \reali;\reali)$,
  $u_o \in \L\infty (\reali^n; \reali)$ and
  $u \in \L\infty (I \times \reali^n; \reali)$. Choose
  $a,b \in \reali$ so that $[a,b] \supseteq J_{u,u_o}$. If
  \begin{displaymath}
    E \in \C2 (I \times \reali^n \times \reali;\reali) \cap \mathscr{F}_c
    \; \mbox{ and } \;
    \mbox{for all } v \in [a,b] \quad
    (t,x) \mapsto \partial^2_{uu} E (t,x,v)
    \in \C1 (I \times \reali^n; \reali) ,
  \end{displaymath}
  then with the notation in Definition~\ref{def:muM} and in
  Definition~\ref{def:muk},
  \begin{equation}
    \label{eq:17}
    \mathcal{M}_u^{tx} (E)
    =
    \dfrac12 \,
    \int_a^b
    \mu_k \left(\partial^2_{uu} E (\cdot, \cdot, k) \right)
    \d{k}
    +
    \mathcal{M}_u^{tx}
    \left(
      \dfrac{\partial_u E (\cdot,\cdot,a) + \partial_u E (\cdot,\cdot,b)}{2}
      \otimes \Id_{\reali} \right)  \,.
  \end{equation}
\end{theorem}
\section{Technical Details}
\label{sec:technical-details}

Throughout, we denote the Lebesgue measure by $\mathcal{L}$ and the
$\sigma$-algebra of Borel sets by $\mathcal{B}$. The characteristic
function of a set $A$ is $\caratt{A}$, so that
$\caratt{A} (x) = 1 \iff x \in A$. Concerning a function $\phi$
compactly supported in $I \times \reali^n$, we mean that the closure
of
$\left\{(t,x) \in I \times \reali^n \colon \phi (t,x) \neq 0\right\}$
is a compact subset of $I \times \reali^n$. In particular,
$\phi (0,x)$ need not vanish.

\subsection{Proofs Related to Section~\ref{sec:time-space-dependent}}
\label{sec:proofs3}

\begin{proofof}{Proposition~\ref{prop:1}}
  Consider first the different items of Proposition~\ref{prop:1}
  separately.

  \paragraph{\ref{item:9}:} Note that~\eqref{eq:75} defines $F$ only
  up to a function of $t$ and $x$. Hence,
  ${F} (t,x,u) - \check{F} (t,x,u) = \Delta (t,x)$ for a
  $\Delta \in \Cc1 (I \times \reali^n;\reali^n)$. Moreover,
  \begin{eqnarray*} %
    {\mathcal{M}}_u^{tx} (E) - \check{\mathcal{M}}_u^{tx} (E) %
    & = %
    & \int_{I} \int_{\reali^n} \left[ (\div_x {F})\left(t, x, u
      (t,x)\right) - (\div_x \check{F})\left(t, x, u (t,x)\right)
      \right] \d{x} \d{t}
    \\
    & = %
    & \int_{I} \int_{\reali^n} \div_x \Delta (t, x) \d{x} \d{t}
    \\
    & = %
    & 0 %
  \end{eqnarray*}
  since for all $t \in I$, the map $x \mapsto \Delta(t,x)$ is
  compactly supported.

  \paragraph{\ref{item:34}:} This proof directly follows
  from~\eqref{eq:13} and~\eqref{eq:40}.

  \paragraph{\ref{item:1}:} The linearity at~\ref{item:1} follows by
  construction. To prove the second part of the statement note that
  $\rest{E}{I \times \reali^n \times J_{u,u_o}} = 0$ and~\eqref{eq:5}
  directly imply that for a.e.~$t \in I$, $x \in \reali^n$,
  $E\left(t, x, u (t,x)\right) = 0$,
  $\partial_u E\left(t, x, u (t,x)\right)=0$,
  $(\partial_t E)\left(t, x, u (t,x)\right)=0$ and
  $E\left(0, x, u_o (x)\right) = 0$. Hence, from
  Definition~\ref{def:muM},
  \begin{displaymath}
    \mathcal{M}_u^{tx} (E)
    =
    \int_{I} \int_{\reali^n}
    (\div_x F)\left(t, x, u (t,x)\right)
    \d{x} \d{t} \,.
  \end{displaymath}
  Call $\bar u \colonequals \inf J_u$. Then, by~\ref{item:9} and
  Definition~\ref{def:EEF}, we can assume that
  $F (t,x,w) = \int_{\bar u}^w \partial_u E (t,x,k) \; \partial_u f
  (t, x,k) \d{k}$ which ensures that for
  a.e.~$(t,x) \in I\times \reali^n$ and for all $w \in J_{u,u_o}$,
  $F(t,x,w)=0$. Hence, for a.e.~$(t,x) \in I\times \reali^n$,
  $(\div_x F)\left(t, x, u (t,x)\right)$ vanishes, completing the
  proof of the claim.
\end{proofof}

\begin{proofof}{Theorem~\ref{thm:continuity}}
  By~\ref{item:9} in Proposition~\ref{prop:1}, it is sufficient to
  consider the flux $F_\nu \colonequals \mathcal{F} (E_\nu,f_\nu)$,
  where we used the notation~\eqref{eq:44}. By Lemma~\ref{lem:treves},
  there exists a compact $K_{tx} \subset I \times \reali^n$ such that
  for all $\nu \in \naturali$, ${\spt}_{tx} E_\nu$, ${\spt}_{tx}E$,
  ${\spt}_{tx}F_\nu$ and ${\spt}_{tx}F$ are all contained in
  $K_{tx}$. Moreover, there exists a compact set $K_u \in \reali$ such
  that for all $\nu \in \naturali$, $u_o (\reali)$,
  $u_o^\nu (\reali)$, $u^\nu (I\times\reali)$ and
  $u_\infty (I\times\reali)$ are all contained in $K_u$.

  Referring to the terms in Definition~\ref{def:muM}, we have the
  uniform convergence on $K_{tx} \times K_u$ of
  $\partial_u E_\nu \underset{\nu\to+\infty}{\longrightarrow}
  \partial_u E$,
  $\div_x f_\nu \underset{\nu\to+\infty}{\longrightarrow} \div_x f$,
  $\div_x F_\nu \underset{\nu\to+\infty}{\longrightarrow} \div_x F$,
  $\partial_t E_\nu \underset{\nu\to+\infty}{\longrightarrow}
  \partial_t E$ and $g \underset{\nu\to+\infty}{\longrightarrow}
  g$. This, combined with the pointwise a.e.~convergence
  $u_o^\nu \underset{\nu\to+\infty}{\longrightarrow} u_o$ and
  $u^\nu \underset{\nu\to+\infty}{\longrightarrow} u_\infty$, allows
  to apply the Dominated Convergence Theorem, completing the proof.
\end{proofof}

\begin{remark}
  \label{rem:NoDer}
  \rm Note for later use that in the proof of the first part
  of~\ref{item:1} in Proposition~\ref{prop:1}, the uniform convergence
  of $E'_\nu$ plays a key role, while its continuity has never been
  exploited.
\end{remark}

\begin{lemma}
  \label{lem:MminusM}
  Let $f \in \mathscr{F}^n$,
  $g \in \C0 (\overline{I}\times \reali^n\times \reali;\reali)$
  satisfying~\eqref{eq:33} and $u_o \in \L\infty (\reali^n, \reali)$.
  Fix $u_1,u_2 \in \L\infty (I \times \reali^n; \reali)$ with
  $J_{u_1,u_o} = J_{u_2,u_o}$ in the sense of~\eqref{eq:5}. Then, for
  any
  $E \in \C2 (I \times \reali^n \times \reali; \reali) \cap
  \mathscr{F}_c$ setting $\psi = \partial_u E$ and using the
  notation~\eqref{eq:32}, the following identity holds:
  \begin{eqnarray*}
    &
    & \mathcal{M}_{u_1}^{tx} (E) - \mathcal{M}_{u_2}^{tx} (E)
    \\
    & =
    & \int_{I} \int_{\reali^n} \int_{\reali}
      v (t,x,\xi)
      \bigl(
      \partial_t \psi (t,x,\xi) \,
      + \partial_u \left( g (t,x,\xi) \, \psi (t,x,\xi) \right)
    \\
    &
    & \qquad \qquad \qquad
      + \div_x \left(\partial_u f (t,x,\xi) \, \psi (t,x,\xi)\right)
      - \partial_u \left(\div_x f (t,x,\xi) \, \psi (t,x,\xi)\right)
      \bigr)
      \d\xi \d{x} \d{t} \,.
  \end{eqnarray*}
\end{lemma}

\begin{proofof}{Lemma~\ref{lem:MminusM}}
  Using Definition~\ref{def:muM}, compute
  \begin{eqnarray}
    \nonumber
    &
    & \mathcal{M}_{u_1}^{tx} (E) - \mathcal{M}_{u_2}^{tx} (E)
    \\
    \label{eq:28}
    &
    & \left.%
      \begin{array}{@{}l@{}}
        \displaystyle
        = - \int_{I} \int_{\reali^n}
        (\partial_u E)\left(t, x, u_1 (t,x)\right)
        \; (\div_x f)\left(t, x, u_1 (t,x)\right)
        \d{x} \d{t}
        \\
        \displaystyle+
        \int_{I} \int_{\reali^n}
        (\partial_u E)\left(t, x, u_2 (t,x)\right)
        \; (\div_x f)\left(t, x, u_2 (t,x)\right)
        \d{x} \d{t}
      \end{array}
      \right\}
    \\
    \label{eq:29}
    &
    & +
      \int_{I} \int_{\reali^n}
      \left(
      (\div_x F)\left(t, x, u_1 (t,x)\right)
      -
      (\div_x F)\left(t, x, u_2 (t,x)\right)
      \right)
      \d{x} \d{t}
    \\
    \label{eq:30}
    &
    & +
      \int_{I} \int_{\reali^n}
      \left(
      (\partial_t E)\left(t, x, u_1 (t,x)\right)
      -
      (\partial_t E)\left(t, x, u_2 (t,x)\right)
      \right)
      \d{x} \d{t}
    \\
    \label{eq:31}
    &
    & \left.%
      \begin{array}{@{}l@{}}
        \displaystyle
        +
        \int_{I} \int_{\reali^n}
        (\partial_u E)\left(t, x, u_1 (t,x)\right) \;
        g \left(t,x,u_1 (t,x)\right)
        \d{x} \d{t}
        \\
        \displaystyle
        -
        \int_{I} \int_{\reali^n}
        (\partial_u E)\left(t, x, u_2 (t,x)\right) \;
        g \left(t,x,u_2 (t,x)\right)
        \d{x} \d{t}
      \end{array}
      \right\}
  \end{eqnarray}
  Using~\eqref{eq:51}, rewrite the terms~\eqref{eq:28}--\eqref{eq:31}
  as follows:
  \begin{align*}
    [\mbox{\eqref{eq:28}}]
    =
    & - \int_{I} \int_{\reali^n} \int_{\reali}
      \partial^2_{uu} E (t,x,\xi) \, \div_x f (t,x,\xi)
 \, v (t,x,\xi)
      \d\xi \d{x} \d{t}
    \\
    & - \int_{I} \int_{\reali^n} \int_{\reali}
      \partial_u E (t,x,\xi) \, \div_x \partial_u f (t,x,\xi)
 \, v (t,x,\xi)
      \d\xi \d{x} \d{t};
    \\
    [\mbox{\eqref{eq:29}}]
    =
    & \int_{I} \int_{\reali^n} \int_{\reali}
      \div_x \partial_u F (t,x,\xi) \, v (t,x,\xi)
      \d\xi \d{x} \d{t}
    \\
    =
    & \int_{I} \int_{\reali^n} \int_{\reali}
      \div_x \left(\partial_u f (t,x,\xi) \, \partial_u E (t,x,\xi)\right)
      \, v (t,x,\xi)
      \d\xi \d{x} \d{t}
    \\
    [\mbox{\eqref{eq:30}}]
    =
    & \int_{I} \int_{\reali^n} \int_{\reali}
      \partial_t \partial_u E (t,x,\xi) \, v (t,x,\xi)
      \d\xi \d{x} \d{t}
    \\
    [\mbox{\eqref{eq:31}}]
    =
    & \int_{I} \int_{\reali^n} \int_{\reali}
      \partial^2_{uu} E (t,x,\xi) \, g (t,x,\xi) \, v (t,x,\xi)
      \d\xi \d{x} \d{t}
    \\
    & + \int_{I} \int_{\reali^n} \int_{\reali}
      \partial_u E (t,x,\xi) \, \partial_u g (t,x,\xi)
      \, v (t,x,\xi)
      \d\xi \d{x} \d{t}
    \\
    =
    & \int_{I} \int_{\reali^n} \int_{\reali}
      \partial_u \left(\partial_u E (t,x,\xi)\right)
      \, g (t,x,\xi) \, v (t,x,\xi)
      \d\xi \d{x} \d{t}
    \\
    & + \int_{I} \int_{\reali^n} \int_{\reali}
      \left(\partial_u E (t,x,\xi) \right)
      \, \partial_u g (t,x,\xi)
      \, v (t,x,\xi)
      \d\xi \d{x} \d{t}
  \end{align*}
  while on the other hand
  \begin{eqnarray*}
    [\mbox{\eqref{eq:28}}] + [\mbox{\eqref{eq:29}}]
    & =
    & - \int_{I} \int_{\reali^n} \int_{\reali} \partial_u
      \left(\partial_u E (t,x,\xi) \right) \, \div_x f (t,x,\xi) \, v
      (t,x,\xi) \d\xi \d{x} \d{t}
    \\
    &
    & +
      \int_{I} \int_{\reali^n} \int_{\reali}
      \partial_u f (t,x,\xi) \, \cdot
      \grad_x \left(\partial_u E (t,x,\xi) \right) \, v (t,x,\xi)
      \d\xi \d{x} \d{t}.
  \end{eqnarray*}
  Combining the expressions obtained for the
  terms~\eqref{eq:28}--\eqref{eq:31}, we have that for
  $\psi = \partial_u E$
  \begin{eqnarray*}
    &
    & \mathcal{M}_{u_1}^{tx} (E)
      -
      \mathcal{M}_{u_2}^{tx} (E)
    \\
    & =
    & \int_{I} \int_{\reali^n} \int_{\reali}
      v (t,x,\xi)
      \bigl(
      \partial_t \psi (t,x,\xi) \,
      + \partial_u \left( g (t,x,\xi) \, \psi (t,x,\xi) \right)
    \\
    &
    & \qquad \qquad \qquad
      + \partial_u f (t,x,\xi) \cdot \grad_x \psi (t,x,\xi)
      - \div_x f (t,x,\xi) \, \partial_u \psi (t,x,\xi)
      \bigr)
      \d\xi \d{x} \d{t}
    \\
    & =
    & \int_{I} \int_{\reali^n} \int_{\reali}
      v (t,x,\xi)
      \bigl(
      \partial_t \psi (t,x,\xi) \,
      + \partial_u \left( g (t,x,\xi) \, \psi (t,x,\xi) \right)
    \\
    &
    & \qquad \qquad \qquad
      + \div_x \left(\partial_u f (t,x,\xi) \, \psi (t,x,\xi)\right)
      - \partial_u \left(\div_x f (t,x,\xi) \, \psi (t,x,\xi)\right)
      \bigr)
      \d\xi \d{x} \d{t}
  \end{eqnarray*}
  completing the proof.
\end{proofof}

\begin{proofof}{Theorem~\ref{thm:Identify}}
  Apply Lemma~\ref{lem:MminusM}, obtaining that for any
  $E \in \C2 (I \times \reali^n \times \reali; \reali) \cap
  \mathscr{F}_c$,
  \begin{equation}
    \label{eq:35}
    \begin{array}{@{}r@{\;}c@{\;}l@{}}
      0
      & =
      & \displaystyle
        \int_{\reali_+} \int_{\reali^n} \int_{\reali}
        v (t,x,\xi)
        \bigl(
        \partial_t \psi (t,x,\xi) \,
        + \partial_u \left( g (t,x,\xi) \, \psi (t,x,\xi) \right)
      \\
      &
      & \qquad \qquad
        \displaystyle
        + \div_x \left(\partial_u f (t,x,\xi) \, \psi (t,x,\xi)\right)
        - \partial_u \left(\div_x f (t,x,\xi) \, \psi (t,x,\xi)\right)
        \bigr)
        \d\xi \d{x} \d{t} \,.
    \end{array}
  \end{equation}
  Given any test function
  $\psi \in \mathcal{D} (I \times \reali^n \times \reali; \reali)$,
  set $E (t,x,u) = \int_0^u \psi (t,x,\xi) \d\xi$ so that
  $E \in \C2 (I \times \reali^n \times \reali; \reali) \cap
  \mathscr{F}_c$ and $\psi = \partial_u E$.

  The validity of~\eqref{eq:35} for any
  $\psi \in \mathcal{D} (I \times \reali^n \times \reali; \reali)$
  shows that $v$ is a weak solution to the following Cauchy problem
  for a linear transport equation
  \begin{equation}
    \label{eq:34}
    \left\{
      \begin{array}{l}
        \partial_t v + \partial_u f \cdot \grad_x v + (g-\div_x f) \, \partial_u v
        = 0
        \\
        v (0,x,\xi) = 0
      \end{array}
    \right.
    \qquad t \in I \,,\, x \in \reali^n \,,\, \xi \in \reali \,.
  \end{equation}
  By~\eqref{eq:33}, we can apply~\cite[Corollary~II.1]{MR1022305} to
  prove that the only solution to~\eqref{eq:34} is the function
  $(t,x,\xi) \mapsto 0$. Hence, by~\eqref{eq:25}, \eqref{eq:51} with
  $\eta (u) = u$ and~\eqref{eq:32}, $u_1 = u_2$ a.e.~in
  $I \times \reali^n$, completing the proof of
  Theorem~\ref{thm:Identify}.
\end{proofof}

\begin{proofof}{Proposition~\ref{prop:affine}}
  With the notation in Proposition~\ref{prop:fEF}, let
  $F \colonequals \mathcal{F} (E,f)$, so that $F \in \mathscr{F}^n_c$
  and $(E,F)$ is an entropy-entropy flux pair as in
  Definition~\ref{def:EEF}.

  \noindent\textbf{Proof of~\ref{item:31}:} Also $(\tilde E, F)$
  is an entropy-entropy flux pair as in Definition~\ref{def:EEF}. By
  direct computations:
  \begin{displaymath}
    \mathcal{M}_u^{tx} (E)
    -
    \mathcal{M}_u^{tx} (\tilde E)
    =
    \int_{I} \int_{\reali^n}
    \partial_t b(t, x)  \d{x} \d{t}
    +\int_{\reali^n} b(0, x)  \d{x}
  \end{displaymath}
  which vanishes since $b$ is compactly supported in
  $I \times \reali^n$.

  \noindent\textbf{Proof of~\ref{item:32}} The map $\tilde F (t,x,u) \colonequals F (t,x,u) + a (t,x) \; f (t,x,u)$
  is such that $\tilde F \in \mathscr{F}^n_c$ and
  $(\tilde E, \tilde F)$ is an entropy-entropy flux pair. Moreover, by
  Definition~\ref{def:muM},
  \begin{eqnarray*}
    &
    & \mathcal{M}_u^{tx} (E) - \mathcal{M}_u^{tx} (\tilde E)
    \\
    & =
    & \int_{\reali_+} \int_{\reali^n}
      u (t,x) \, \partial_t a (t,x)
      \d{x} \d{t}
      +
      \int_{\reali_+} \int_{\reali^n}
      f \left(t,x,u (t,x) \right) \cdot
      \grad_x a (t,x)
      \d{x} \d{t}
    \\
    &
    & +
      \int_{\reali_+} \int_{\reali^n}
      g\left(t,x,u (t,x)\right) a (t,x)
      \d{x} \d{t}
      +
      \int_{\reali^n}
      u_o (x)\, a (0,x) \d{x}
    \\
    &
    & +
      \int_{\reali_+} \int_{\reali^n}
      \partial_t b (t,x)
      \d{x} \d{t}
      +
      \int_{\reali^n}
      b (0,x) \d{x}
    \\
    & =
    & 0
  \end{eqnarray*}
  where we used $a$ as a test function in Definition~\ref{def:ws} and
  the fact that $b$ is compactly supported in $I \times \reali^n$.
\end{proofof}

\begin{proofof}{Theorem~\ref{thm:formulaC2exProp}}
  Consider first the case
  $E \in \C3( I \times \reali^n \times \reali; \reali) \cap
  \mathscr{F}_c$. By Fubini Theorem we have
  \begin{eqnarray*}
    &
    & \int_a^b \mu_k
      \left(\partial^2_{uu} E (\cdot, \cdot, k)\right) \d{k}
    \\
    & =
    &  \int_a^b\!\!\!
      \int_{I}\! \int_{\reali^n}
      \left(
      \modulo{u (t,x) - k} \, \partial_t \partial^2_{uu} E(t,x,k)
      +
      \Phi\left(t,x, u (t,x),k\right) \cdot \grad_x \partial^2_{uu} E(t,x,k)
      \right) \d{x} \d{t} \d{k}
    \\
    &
    & -
      \int_a^b
      \int_{I} \int_{\reali^n}
      \partial^2_{uu} E (t,x,k) \,
      \sgn\left(u (t,x) - k\right) \; \div_x f (t,x,k)
      \d{x} \d{t} \d{k}
    \\
    &
    & +
      \int_a^b
      \int_{I} \int_{\reali^n}
      \partial^2_{uu} E (t,x,k) \,
      \sgn\left(u (t,x) - k\right) \; g \left(t,x,u(t,x)\right)
      \d{x} \d{t} \d{k}
    \\
    &
    & +
      \int_a^b
      \int_{\reali^n}
      \partial^2_{uu} E (0,x,k) \,\modulo{u_o(x) - k} \d{x}
      \d{k}
    \\
    & =
    & \int_{I} \int_{\reali^n}
      \int_a^b
      \modulo{u (t,x) - k} \,  \partial^2_{uu} \partial_t E(t,x,k)
      \d{k} \d{x} \d{t}
    \\
    &
    & + \int_{I} \int_{\reali^n}
      \int_a^b
      \Phi\left(t,x, u (t,x),k\right) \cdot \partial^2_{uu} \grad_x E(t,x,k)
      \d{k} \d{x} \d{t}
    \\
    &
    & -
      \int_{I} \int_{\reali^n}
      \int_a^b \partial^2_{uu}E (t,x,k) \; \sgn\left(u (t,x) - k\right) \;
      \div_x f (t,x,k)
      \d{k} \d{x} \d{t}
    \\
    &
    & +
      \int_{I} \int_{\reali^n}
      \int_a^b \partial^2_{uu}E (t,x,k)
      \;  \sgn\left(u (t,x) - k\right) \; g \left(t,x,u(t,x)\right)
      \d{k} \d{x} \d{t}
    \\
    &
    & +
      \int_{\reali^n}
      \int_a^b \partial^2_{uu}E (0,x,k)
      \; \modulo{u_o(x) - k}
      \d{k} \d{x}
    \\
    & =
    & \int_{I} \int_{\reali^n}
      A_1 (t,x) \d{x} \d{t}
      +
      \int_{I} \int_{\reali^n}
      A_2 (t,x) \d{x} \d{t}
    \\
    &
    & -
      \int_{I} \int_{\reali^n}
      A_3 (t,x) \d{x} \d{t}
      +
      \int_{I} \int_{\reali^n}
      A_4 (t,x) \d{x} \d{t}
      +
      \int_{\reali^n}
      A_5 (x) \d{x}
  \end{eqnarray*}
  and consider the terms above separately. Using Schwarz Lemma and
  recalling that $u (t,x) \in [a,b]$ for
  a.e.~$(t,x) \in I\times\reali^n$,
  \begin{eqnarray*}
    A_1 (t,x)
    & =
    & \int_a^{u (t,x)}
      \left(u (t,x) - k\right) \,  \partial^2_{uu} \partial_t E(t,x,k)
      \d{k}
      -
      \int_{u (t,x)}^b
      \left(u (t,x) - k\right) \,  \partial^2_{uu} \partial_t E(t,x,k)
      \d{k}
    \\
    & =
    & - \left(u (t,x) - a\right) \, \partial_u \partial_t E(t,x,a)
      + \int_a^{u (t,x)} \partial_u \partial_t E(t,x,k) \d{k}
    \\
    &
    & - \left(u (t,x) - b\right)\, \partial_u \partial_t E(t,x,b)
      -\int_{u (t,x)} ^ b \partial_u \partial_t E(t,x,k) \d{k}
    \\
    & =
    & - \left(u (t,x) - a\right) \, \partial_u \partial_t E(t,x,a)
      +  \partial_t E\left(t,x,u (t,x)\right)
      -  \partial_t E(t,x,a)
    \\
    &
    & - \left(u (t,x) - b\right)\, \partial_u \partial_t E(t,x,b)
      + \partial_t E\left(t,x,u (t,x)\right)
      -  \partial_t E(t,x,b)
    \\
    & =
    & 2 \, \partial_t E\left(t,x,u (t,x)\right)
      - \left(u (t,x) - a\right) \, \partial_u \partial_t E(t,x,a)
      -  \partial_t E(t,x,a)
    \\
    &
    & - \left(u (t,x) - b\right)\, \partial_u \partial_t E(t,x,b)
      -  \partial_t E(t,x,b) \,.
  \end{eqnarray*}
  Repeat now similar computations:
  \begin{eqnarray*}
    A_2 (t,x)
    & =
    & \int_a^{u (t,x)} \partial^2_{uu} \grad_x E (t,x,k) \cdot
      \Phi\left(t,x, u (t,x),k\right) \d{k}
    \\
    &
    & +
      \int_{u (t,x)}^b \partial^2_{uu} \grad_x E (t,x,k) \cdot
      \Phi\left(t,x, u (t,x),k\right) \d{k}
    \\
    & =
    & \int_a^{u (t,x)} \partial^2_{uu} \grad_x E (t,x,k) \cdot
      \left(f\left(t,x, u (t,x)\right) - f (t,x,k)\right) \d{k}
    \\
    &
    & \quad
      -
      \int_{u (t,x)}^b \partial^2_{uu} \grad_x E (t,x,k) \cdot
      \left(f\left(t,x, u (t,x)\right) - f (t,x,k)\right) \d{k}
    \\
    & =
    & - \partial_u \grad_x E (t,x,a) \cdot
      f\left(t,x, u (t,x)\right)
      +
      \partial_u \grad_x E (t,x,a) \cdot f (t,x,a)
    \\
    &
    & + \int_a^{u (t,x)} \partial_u \grad_x E (t,x,k) \cdot
      \partial_u f (t,x,k) \d{k}
    \\
    &
    & - \partial_u \grad_x E (t,x,b) \cdot
      f\left(t,x, u (t,x)\right)
      +
      \partial_u \grad_x E (t,x,b) \cdot f (t,x,b)
    \\
    &
    & - \int_{u (t,x)}^b \partial_u \grad_x E (t,x,k) \cdot
      \partial_u f (t,x,k) \d{k} \,.
  \end{eqnarray*}
  By~\eqref{eq:75} and using Schwarz Lemma
  \begin{eqnarray*}
    A_3 (t,x)
    & =
    & \int_a^{u (t,x)} \partial^2_{uu} E (t,x,k) \,
      \sgn\left(u (t,x) - k\right) \,
      \div_x f (t,x,k) \d{k}
    \\
    &
    & \quad
      +
      \int_{u (t,x)}^b \partial^2_{uu} E (t,x,k) \,
      \sgn\left(u (t,x) - k\right) \,
      \div_x f (t,x,k) \d{k}
    \\
    & =
    & \int_a^{u (t,x)} \partial^2_{uu} E (t,x,k) \,
      \div_x f (t,x,k) \d{k}
      -
      \int_{u (t,x)}^b \partial^2_{uu} E (t,x,k) \,
      \div_x f (t,x,k) \d{k}
    \\
    & =
    & \partial_u E\left(t, x, u (t,x)\right) \,
      \div_x f\left(t,x, u (t,x)\right)
      - \partial_u E (t,x,a) \, \div_x f (t,x,a)
    \\
    &
    & \quad -
      \int_a^{u (t,x)} \partial_u E (t,x,k)\, \partial_u \div_x f (t,x,k) \d{k}
    \\
    &
    & \quad +
      \partial_u E\left(t, x, u (t,x)\right) \,
      \div_x f\left(t,x, u (t,x)\right)
      - \partial_u E (t,x,b) \, \div_x f (t,x,b)
    \\
    &
    & \quad +
      \int_{u (t,x)}^b \partial_u E (t,x,k)\; \partial_u \div_x f (t,x,k) \d{k}\,.
  \end{eqnarray*}
  Combining the expressions obtained for $A_2$ and $A_3$,
  by~\eqref{eq:75} we obtain:
  \begin{eqnarray*}
    &
    & A_2 (t,x) - A_3 (t,x)
    \\
    & =
    & - \partial_u \grad_x E (t,x,a) \cdot
      f\left(t,x, u (t,x)\right)
      +
      \partial_u \grad_x E (t,x,a) \cdot f (t,x,a)
    \\
    &
    & - \partial_u \grad_x E (t,x,b) \cdot
      f\left(t,x, u (t,x)\right)
      +
      \partial_u \grad_x E (t,x,b) \cdot f (t,x,b)
    \\
    &
    & - \partial_u E\left(t, x, u (t,x)\right) \,
      \div_x f\left(t,x, u (t,x)\right)
      + \partial_u E (t,x,a) \, \div_x f (t,x,a)
    \\
    &
    & - \partial_u E\left(t, x, u (t,x)\right) \,
      \div_x f\left(t,x, u (t,x)\right)
      + \partial_u E (t,x,b) \, \div_x f (t,x,b)
    \\
    &
    & + \int_a^{u (t,x)} \partial_u \div_x F (t,x,k) \d{k}
      - \int_{u (t,x)}^b \partial_u \div_x F (t,x,k) \d{k}
    \\
    & =
    & - 2\, \partial_u E\left(t, x, u (t,x)\right) \,
      \div_x f\left(t,x, u (t,x)\right)
      + 2 \, \div_x F\left(t,x,u (t,x)\right)
    \\
    &
    & - \left(
      \grad_x \partial_u E (t,x,a) +\grad_x \partial_u E (t,x,b)
      \right) f\left(t,x,u (t,x)\right)
    \\
    &
    & + \partial_u \grad_x E (t,x,a) \cdot f (t,x,a)
      + \partial_u E (t,x,a) \, \div_x f (t,x,a)
    \\
    &
    & + \partial_u \grad_x E (t,x,b) \cdot f (t,x,b)
      + \partial_u E (t,x,b) \, \div_x f (t,x,b)
    \\
    &
    & - \div_x F(t,x,a) - \div_x F(t,x,b)
  \end{eqnarray*}
  Similarly,
  \begin{eqnarray*}
    A_4 (t,x)
    & =
    &
      \left(
      \int_a^{u (t,x)} \partial^2_{uu} E (t,x,k) \d{k}
      -
      \int_{u (t,x)}^b \partial^2_{uu} E (t,x,k) \d{k}
      \right)
      g \left(t,x,u(t,x)\right)
    \\
    & =
    & \left(2 \, \partial_u E\left(t,x,u (t,x)\right) - \partial_u E (t,x,a)
      - \partial_u E (t,x,b)\right)  \, g \left(t,x,u(t,x)\right) \,.
  \end{eqnarray*}
  Use~\eqref{eq:6} and recall that $u (t,x) \in [a,b]$ for
  a.e.~$(t,x)$ in $I \times \reali^n$.
  \begin{eqnarray*}
    A_5 (x)
    & =
    & \int_a^{u_o (x)} \partial^2_{uu} E (0,x,k) \; \modulo{u_o(x) - k} \d{k}
      +
      \int_{u_o (x)}^b \partial^2_{uu} E (0,x,k) \; \modulo{u_o(x) - k} \d{k}
    \\
    & =
    &  \int_a^{u_o (x)} \partial^2_{uu} E (0,x,k) \; \left(u_o(x) - k\right) \d{k}
      -
      \int_{u_o (x)}^b \partial^2_{uu} E (0,x,k) \; \left(u_o(x) - k\right) \d{k}
    \\
    & =
    & -
      \partial_u E (0,x,a) \left(u_o(x) - a\right)
      +
      \int_a^{u_o(x)} \partial_u E (0,x,k) \d{k}
    \\
    &
    & \qquad
      -
      \int_{u_o(x)}^b \partial_u E (0,x,k) \d{k}
      -
      \partial_u E (0,x,b) \left(u_o(x) - b\right)
    \\
    & =
    & 2 E\left(0,x,u_o(x)\right)
      -
      u_o (x) \left(\partial_u E (0,x,a) + \partial_u E (0,x,b)\right)
    \\
    &
    & \qquad
      +
      a \, \partial_u E (0,x,a) + b\, \partial_u E (0,x,b)
      - E (0,x,a) - E (0,x,b) \,.
  \end{eqnarray*}
  Adding up all the expressions obtained, by Definition~\ref{def:muM},
  we have:
  \begin{flalign*}
    & \int_a^b \mu_k \left(\partial^2_{uu} E (\cdot, \cdot, k) \right) \d{k}
    \\
    = & 2\, \mathcal{M}_u^{tx} (E)
    \\
    & \left.
      \begin{array}{@{}l}
        \displaystyle
        - \Bigl[ \int_{I} \int_{\reali^n} \left( \partial_u E
        (t,x,a) + \partial_u E (t,x,b) \right)
        g\left(t,x,u (t,x)\right)  \d{x} \d{t}
        \\
        \displaystyle
        \quad + \int_{\reali^n} \left( \partial_u E(0,x,a) + \partial_u
        E(0,x,b) \right) u_o (x) \d{x}
        \\
        \displaystyle
        \quad + \int_{I} \int_{\reali^n} \left( \partial_t
        \partial_u E (t,x,a) + \partial_t \partial_u E (t,x,b) \right) u
        (t,x) \d{x} \d{t}
        \\
        \displaystyle
        \quad + \int_{I} \int_{\reali^n} \left( \grad_x
        \partial_u E (t,x,a) + \grad_x \partial_u E (t,x,b) \right)
        f\left(t, x, u(t,x)\right) \d{x} \d{t} \Bigr]
      \end{array}
      \right\} B_1
    \\
    & \left.
      \begin{array}{@{}l}
        \displaystyle
        + \Bigl[ \int_{I} \int_{\reali^n} \partial_t \left( a
        \,\partial_u E (t,x,a) - E (t,x,a) + b \,\partial_u E (t,x,b) -
        E (t,x,b) \right) \d{x}\d{t}
        \\
        \displaystyle
        \quad + \int_{\reali^n} \left( a \,\partial_u E (0,x,a) - E
        (0,x,a) + b \,\partial_u E (0,x,b) - E (0,x,b) \right) \d{x} \Bigr]
      \end{array}
      \right\} B_2
    \\
    & \left.
      \begin{array}{@{}l@{}}
        \displaystyle
        + \int_{I} \int_{\reali^n} \Bigl[ \partial_u
        \grad_x E (t,x,a) \cdot f (t,x,a) + \partial_u E (t,x,a) \, \div_x
        f (t,x,a) - \div_x F (t,x,a)
        \\
        \displaystyle
        \qquad + \partial_u \grad_x E (t,x,b) \cdot f
        (t,x,b) + \partial_u E (t,x,b) \, \div_x f (t,x,b) - \div_x F
        (t,x,b) \Bigr] \! \d{x} \d{t}
      \end{array}
      \right\} B_3
    \\
    = & 2\, \mathcal{M}_u^{tx} (\phi) - B_1 + B_2 + B_3 \,.
  \end{flalign*}
  To compute now the terms $B_1,B_2,B_3$, introduce
  \begin{eqnarray*}
    \psi_1 (t,x)
    & \colonequals
    & \partial_u E (t,x,a) + \partial_u E (t,x,b)
    \\
    \psi_2 (t,x)
    & \colonequals
    & a \,\partial_u E (t,x,a) - E (t,x,a) + b
      \,\partial_u E (t,x,b) - E (t,x,b)
    \\
    \psi_3 (t,x)
    & \colonequals
    & \partial_u E (t,x,a) \, f (t,x,a) - F (t,x,a)
      +
      \partial_u E (t,x,b) \, f (t,x,b) - F (t,x,b) \,.
  \end{eqnarray*}
  They are all of class $\C1$ and compactly supported since
  $E \in \mathscr{F}_c$ and $F \in \mathscr{F}^n_c$. So, they can be
  used as test function. Hence,
  \begin{eqnarray*}
    B_1
    & =
    & \int_{I} \int_{\reali^n}
      \left(
      u (t,x) \, \partial_t \psi_1 (t,x)
      +
      f\left(t, x, u (t,x)\right) \cdot \grad_x \psi_1 (t,x)
      \right) \d{x} \d{t}
    \\
    &
    & +
      \int_{I} \int_{\reali^n}
      g \left(t,x,u (t,x)\right) \; \psi_1 (t,x)
      \d{x} \d{t}
      +\int_{\reali^n} u_o (x) \, \psi_1 (0,x) \d{x}
    \\
    & =
    & \mathcal{M}_u^{tx} (\psi_1 \otimes \Id_{\reali}) \,.
    \\
    B_2
    & =
    & \int_{I} \int_{\reali^n} \partial_t \psi_2 (t,x) \d{x}\d{t}
      + \int_{\reali^n} \psi_2 (0,x) \d{x}
      = 0 \,.
    \\
    B_3
    & =
    & \int_{I} \int_{\reali^n} \div_x \psi_3 (t,x) \d{x} \d{t}
      = 0 \,.
  \end{eqnarray*}
  The proof is completed if
  $E \in \C3 (I \times \reali^n \times \reali; \reali)$. The general
  case follows by the density proved in Corollary~\ref{cor:density},
  thanks to Theorem~\ref{thm:continuity} and to~\ref{item:10} in
  Proposition~\ref{prop:2}.
\end{proofof}

Above, we refer to~\ref{item:10} in Proposition~\ref{prop:2} which is
proved in the next section. Its proof is however entirely independent
from the results in the present section.

\subsection{Proofs Related to \S~\ref{subsec:first-setting}}
\label{sec:proofs2}

\begin{proofof}{Proposition~\ref{prop:1notx}}
  Items~\ref{item:9notx}, \ref{item:33} and~\ref{item:1notx} directly
  follow respectively from~\ref{item:9}, \ref{item:34}
  and~\ref{item:1} by Lemma~\ref{lem:extension}.

  We prove now~\ref{item:4notx}. For any $E \in \C1 (\reali;\reali)$,
  $\mathcal{M}_u (E)$ is linear by construction and its continuity
  follows noting that for any compact
  $K_{tx} \subset I \times \reali^n$ there is a constant $C$ depending
  on $f,g,u_o,u,E,F$ such that
  \begin{displaymath}
    \modulo{\left(\mathcal{M}_u (E)\right) (\phi)}
    \leq
    C \left(\norma{\phi}_{\L\infty (I\times\reali^n;\reali)} +
      \norma{
        \big[\partial_t\phi \quad \grad_x \phi\big]
      }_{\L\infty(I\times\reali^n;\reali^{n+1})}
    \right)
  \end{displaymath}
  for any $\phi \in \C1 (\reali \times \reali^n;\reali)$ with
  $\spt\phi \subseteq K_{tx}$. This also proves that the distribution
  $\mathcal{M}_u (E)$ is at most of order $1$.
\end{proofof}

The next proof shows that Theorem~\ref{thm:continuity} implies
Theorem~\ref{thm:continuitynotx}. More precisely, it also allows
to prove a sort of uniform convergence, in the sense that also $\phi$
might be allowed to vary, proving that
$\left(\mathcal{M}_\nu (E_\nu)\right) (\phi_\nu)
\underset{\nu\to+\infty}{\longrightarrow} \left(\mathcal{M}_{u_\infty}
  (E)\right) (\phi)$.

\begin{proofof}{Theorem~\ref{thm:continuitynotx}}
  The proof consists in showing that for any
  $\phi \in \mathcal{D} (I\times\reali^n; \reali)$, it holds that
  $\left(\mathcal{M}_\nu (E_\nu)\right) (\phi)
  \underset{\nu\to+\infty}{\longrightarrow}
  \left(\mathcal{M}_{u_\infty} (E)\right) (\phi)$. Then,
  \begin{displaymath}
    \begin{array}{rcl@{\qquad\qquad}r}
      \left(\mathcal{M}_\nu (E_\nu)\right) (\phi)
      & =
      & \mathcal{M}^{tx}_\nu (\phi \otimes E_\nu)
      & [\mbox{By Lemma~\ref{lem:extension}}]
      \\
      & \underset{\nu\to+\infty}{\longrightarrow}
      & \mathcal{M}^{tx}_{u_\infty} (\phi \otimes E)
      & [\mbox{By Theorem~\ref{thm:continuity}}]
      \\
      & =
      & \left(\mathcal{M}_{u_\infty} (E)\right) (\phi)
      & [\mbox{By Lemma~\ref{lem:extension}}]
    \end{array}
  \end{displaymath}
  since under the present assumptions, we have the convergence
  $\phi \otimes E_\nu \underset{\nu\to+\infty}{\longrightarrow} \phi
  \otimes E$ in $\mathscr{F}_c$, thanks to Lemma~\ref{lem:treves}.
\end{proofof}

Inspired by~\cite[Chapter~1, \S~1]{MR2064166}, introduce the function
\begin{equation}
  \label{eq:25}
  \begin{array}{ccccc}
    \chi
    & \colon
    & \reali^3
    & \to
    & \reali
    \\
    &
    & (\xi,a,b)
    & \mapsto
    & \left\{
      \begin{array}{ll}
        1
        & b \leq \xi < a,
        \\
        -1
        & a \leq \xi < b,
        \\
        0
        & \mbox{ otherwise.}
      \end{array}
      \right.
  \end{array}
\end{equation}
Note that $\chi$ is chosen so that for all
$\eta \in \AC (\reali; \reali)$, we have
\begin{equation}
  \label{eq:51}
  \forall \, a,b \in \reali \qquad
  \eta (a) - \eta (b)
  =
  \int_{\reali} \eta' (\xi) \; \chi (\xi,a,b) \, \d\xi \,.
\end{equation}
Given arbitrary $u_1,u_2 \in \L\infty (I \times \reali^n; \reali)$, it
is of help to introduce also the function
\begin{equation}
  \label{eq:32}
  \begin{array}{ccccc}
    v
    & \colon
    & I \times \reali^n \times \reali
    & \to
    & \reali
    \\
    &
    & (t,x,\xi)
    & \mapsto
    & \chi\left(\xi,u_1 (t,x), u_2 (t,x)\right) \,.
  \end{array}
\end{equation}

\begin{proofof}{Theorem~\ref{thm:Identifynotx}}
  By Lemma~\ref{lem:extension}, for any $E \in \C1 (\reali;\reali)$
  and any $\phi \in \Cc1 (I\times\reali^N;\reali)$, we have
  \begin{displaymath}
    \mathcal{M}^{tx}_{u_1} (\phi \otimes E)
    =
    \left(\mathcal{M}_{u_1} (E)\right) (\phi)
    =
    \left(\mathcal{M}_{u_2} (E)\right) (\phi)
    =
    \mathcal{M}^{tx}_{u_2} (\phi \otimes E)\,.
  \end{displaymath}
  By the linearity proved in~\ref{item:1} in Proposition~\ref{prop:1},
  $\mathcal{M}^{tx}_{u_1}$ coincides with $\mathcal{M}^{tx}_{u_2}$ on
  the algebraic tensor product
  $\Cc1 (I\times\reali^N; \reali) \otimes \C1 (\reali;\reali)$. By the
  continuity proved in Theorem~\ref{thm:continuity} and the
  density of this the algebraic tensor product in $\mathscr{F}_c$
  proved in Corollary~\ref{cor:density}, we have that
  $\mathcal{M}^{tx}_{u_1}$ coincides with $\mathcal{M}^{tx}_{u_2}$ on
  all $\mathscr{F}_c$. Theorem~\ref{thm:Identify} allows to
  conclude the proof.
\end{proofof}

\begin{proofof}{Proposition~\ref{prop:2}}
  In item~\ref{item:10} the linearity is immediate. Fix compact sets
  $K_t \subset I$, $K_x\subset \reali^n$ and $K_u \in \reali$. Choose
  a test function $\phi \in \mathcal{D} (I\times\reali^N; \reali)$
  with $\spt\phi \subseteq K_t\times K_x$ and use
  Definition~\ref{def:muk} to obtain, for all $k \in K_u$
  \begin{displaymath}
    \modulo{\mu_k (\phi)}
    \leq
    C \left(
      \norma{\phi}_{\C0 (K_{tx};\reali)}
      +
      \norma{\partial_t \phi}_{\C0 (K_{tx};\reali)}
      +
      \norma{\grad_x\phi}_{\C0 (K_{tx};\reali^n)}
    \right)
  \end{displaymath}
  where
  \begin{eqnarray*}
    C
    & \colonequals
    & \mathcal{L} (K_t\times K_x)
      \left(
      \norma{u}_{\L\infty (I\times\reali^N;\reali)}
      +
      \mathcal{L} (K_u)
      +
      \norma{f}_{\L\infty (K_t\times K_x\times J_u);\reali}
      +
      \norma{f}_{\L\infty (K_t\times K_x\times K_u);\reali)}
      \right.
    \\
    &
    & \qquad\qquad\qquad\qquad \left.
      +
      \norma{g}_{\L\infty (K_t\times K_x\times J_u);\reali}
      +
      \norma{\div_x f}_{\L\infty (K_t\times K_x\times K_u);\reali)}
      \right)
    \\
    &
    & +
      \mathcal{L} (K_x)
      \left(\norma{u_o}_{\L\infty (\reali^n; \reali)} + \mathcal{L} (K_u)\right)
  \end{eqnarray*}
  showing that $\mu_k$ is a first order distribution.

  To prove~\ref{item:11}, fix $k_1,k_2 \in K_u$ with $k_1 < k_2$,
  $\phi \in \mathcal{D} (I \times \reali^n;\reali)$ and using also the
  Lipschitz continuity of $\Phi$ in~\eqref{eq:12},
  see~\cite[Lemma~3]{Kruzkov1970}, compute:
  \begin{eqnarray*}
    &
    & \modulo{\mu_{k_1} (\phi) - \mu_{k_2} (\phi)}
    \\
    & \leq
    & \int_{I} \int_{\reali^n}
      \modulo{k_1 - k_2} \, \modulo{\partial_t \phi (t,x)}
      \d{x} \d{t}
    \\
    &
    & +
      \int_{I} \int_{\reali^n}
      \modulo{
      \Phi\left(t, x, u (t,x),k_1\right)-
      \Phi\left(t, x, u (t,x),k_2\right)} \,
      \modulo{\div_x \phi (t,x)}
      \d{x} \d{t}
    \\
    &
    & +
      \int_{I} \int_{\reali^n}
      \modulo{
      \sgn\left(u (t,x) - k_1\right)
      -
      \sgn\left(u (t,x) - k_2\right)
      }
      \modulo{\div_x f (t,x,k_2)}
      \modulo{\phi (t,x)}
      \d{x} \d{t}
    \\
    &
    & +
      \int_{I} \int_{\reali^n}
      \modulo{\sgn\left(u (t,x) - k_1\right)}
      \modulo{
      \div_x f (t,x,k_1)
      -
      \div_x f (t,x,k_2)
      }
      \modulo{\phi (t,x)}
      \d{x} \d{t}
    \\
    &
    & +
      \int_{I} \int_{\reali^n}
      \modulo{\sgn\left(u (t,x) - k_1\right)-\sgn\left(u (t,x) - k_2\right)}
      \, \modulo{g \left(t,x,u (t,x)\right)}
      \, \modulo{\phi (t,x)}
      \d{x} \d{t}
    \\
    &
    & + \int_{\reali^n} \modulo{k_1 - k_2} \, \modulo{\phi(0,x)} \d{x}
    \\
    & \leq
    & \modulo{k_1 - k_2} \int_{I} \int_{\reali^n}
      \modulo{\partial_t \phi (t,x)}
      \d{x} \d{t}
    \\
    &
    & +
      \modulo{k_1 - k_2}
      \int_{I} \int_{\reali^n}
      \norma{\partial_u f}_{\L\infty (\spt\phi \times J_{u,u_o};\reali^n)} \,
      \modulo{\div_x \phi (t,x)}
      \d{x} \d{t}
    \\
    &
    & +
      \int_{I} \int_{\reali^n}
      \left(
      2\caratt{k_1 < u (t,x)<k_2} (t,x)
      + \caratt{k_1 = u (t,x)} (t,x)
      + \caratt{u (t,x)=k_2} (t,x)
      \right)
    \\
    &
    & \qquad \times
      \norma{\div_x f}_{\L\infty (\spt\phi\times K_u;\reali)}
      \modulo{\phi (t,x)}
      \d{x} \d{t}
    \\
    &
    & +
      \modulo{k_1-k_2} \int_{I} \int_{\reali^n}
      \norma{\partial_u \div_x f}_{\L\infty (\spt \phi \times K_u;\reali)}
      \, \modulo{\phi (t,x)}
      \d{x} \d{t}
    \\
    &
    & +
      \int_{I} \int_{\reali^n}
      \left(
      2\caratt{k_1 < u (t,x)<k_2} (t,x)
      + \caratt{k_1 = u (t,x)} (t,x)
      + \caratt{u (t,x)=k_2} (t,x)
      \right)
    \\
    &
    & \qquad \times
      \norma{g}_{\L\infty (\spt\phi\times J_{u,u_o};\reali)}
      \modulo{\phi (t,x)}
      \d{x} \d{t}
    \\
    &
    & + \modulo{k_1 - k_2} \int_{\reali^n} \modulo{\phi(0,x)} \d{x} \,.
  \end{eqnarray*}
  The above computations show that for any
  $\phi \in \mathcal{D} (I \times \reali^n; \reali)$, with
  $\spt\phi \subseteq K_{tx}$, there exists a constant $C$ depending
  on $f,g,u$ such that, denoting
  $K_+ = K_{tx} \cap (I \times \reali^n)$,
  \begin{eqnarray*}
    &
    & \modulo{\mu_{k_1} (\phi) - \mu_{k_2} (\phi)}
    \\
    & \leq
    &
      C \, \norma{\phi}_{\C0 (I \times \reali^n;\reali)}
      \iint_{K_+}
      \left(
      2\caratt{k_1 < u (t,x)<k_2} (t,x)
      + \caratt{k_1 = u (t,x)} (t,x)
      + \caratt{u (t,x)=k_2} (t,x)
      \right)
      \d{x} \d{t}
    \\
    &
    & +
      C \, \mathcal{L} (K_{tx}) \, \norma{\phi}_{\C1 (I \times \reali^n;\reali)}
      \, \modulo{k_1-k_2} \,.
  \end{eqnarray*}
  For any increasing sequence $k_i$ in $K_u$, we then have
  \begin{displaymath}
    \sum_i \modulo{\mu_{k_{i+1}} (\phi) - \mu_{k_i} (\phi)}
    \leq
    C \, \mathcal{L} (K_{tx})
    \left(
      2\, \norma{\phi}_{\C0 (I \times \reali^n;\reali)}
      +
      \mathcal{L} (K_u) \, \norma{\phi}_{\C1 (I \times \reali^n;\reali)}
    \right)\,
  \end{displaymath}
  completing the proof of~\ref{item:11}.

  To prove~\ref{item:27}, note first that the choice~\eqref{eq:6}
  ensures that for all $z \in \reali$, the $\sgn$ function satisfies
  $\sgn z = \frac12 \left(\sgn (z-) + \sgn (z+)\right)$. Then, thanks
  to $\phi$ being bounded and compactly supported and thanks to
  $u,u_o$ being of class $\L\infty$, the Dominated Convergence Theorem
  can be applied to compute the limits in~\ref{item:27} of
  $\mu_k (\phi)$ as written in Definition~\ref{def:muk}. This
  concludes the proof.
\end{proofof}

\begin{proofof}{Proposition~\ref{prop:affinenotx}}
  Let $\phi \in \Cc1 (I\times\reali^N;\reali)$. Then, to
  prove~\ref{item:30}, recall Remark~\ref{rem:confusion} and observe that
  \begin{flalign*}
    \left(\mathcal{M}_u (\tilde E)\right) (\phi)
    & = \mathcal{M}^{tx}_u (\phi \otimes \tilde E)
    & [\mbox{By Lemma~\ref{lem:extension}}]
    \\
    & = \mathcal{M}^{tx}_u (\phi \otimes E) + \mathcal{M}^{tx}_u (\phi \otimes b)
    &[\mbox{By~\ref{item:1} in Proposition~\ref{prop:1}}]
    \\
    & = \mathcal{M}^{tx}_u (\phi \otimes E)
    & [\mbox{By~\ref{item:31} in Proposition~\ref{prop:affine}}]
    \\
    & = \left(\mathcal{M}_u (E)\right) (\phi) \,.
    & [\mbox{By Lemma~\ref{lem:extension}}]
  \end{flalign*}

  To prove~\ref{item:29}, follow the same reasoning:
  \begin{flalign*}
    \left(\mathcal{M}_u (\tilde E)\right) (\phi)
    & = \mathcal{M}^{tx}_u (\phi \otimes \tilde E)
    & [\mbox{By Lemma~\ref{lem:extension}}]
    \\
    & = \mathcal{M}^{tx}_u (\phi \otimes E)
      + \mathcal{M}^{tx}_u \left(\phi \otimes (a\Id_{\reali}+ b)\right)
    &[\mbox{By~\ref{item:1} in Proposition~\ref{prop:1}}]
    \\
    & = \mathcal{M}^{tx}_u (\phi \otimes E)
    & [\mbox{By~\ref{item:32} in Proposition~\ref{prop:affine}}]
    \\
    & = \left(\mathcal{M}_u (E)\right) (\phi) \,,
    & [\mbox{By Lemma~\ref{lem:extension}}]
  \end{flalign*}
  completing the proof.
\end{proofof}

\begin{lemma}
  \label{lem:formulanotx}
  Let $f \in \mathscr{F}^n$,
  $g \in \C0 (\overline{I} \times \reali^n\times \reali;\reali)$,
  $u_o \in \L\infty (\reali^n; \reali)$ and
  $u \in \L\infty (I \times \reali^n; \reali)$. Denote
  $a,b \in \reali$ so that $[a,b] = J_{u,u_o}$. If
  $E \in \C2 (\reali;\reali)$, then with the notation in
  Definition~\ref{def:muMnotx} and in Definition~\ref{def:muk}, for
  all $\phi \in \Cc1 (I\times\reali^n;\reali)$
  \begin{equation}
    \label{eq:19}
    \left(\mathcal{M}_u (E)\right) (\phi)
    =
    \dfrac12 \,
    \int_a^b E'' (k) \, \mu_k (\phi)  \d{k}
    +
    \dfrac{E' (a) + E' (b)}{2}
    \left(\mathcal{M}_u (\Id_{\reali})\right) (\phi)  \,.
  \end{equation}
\end{lemma}

The proof is a direct consequence of Lemma~\ref{lem:extension} and
Theorem~\ref{thm:formulaC2exProp}, using $\phi \otimes E$
in~\eqref{eq:17}.

\begin{proofof}{Theorem~\ref{thm:formulaC2}}
  If~\ref{item:21} holds, then Remark~\ref{rem:identity} show
  equality~\eqref{eq:2}. On the other hand, if~\ref{item:22} holds,
  then the latter term in~\eqref{eq:19} vanishes, choose $E (u) = u$,
  proving~\ref{item:21} by Remark~\ref{rem:identity}.
\end{proofof}

\subsection{Proofs Related to \S~\ref{subsec:case-class-entr}}
\label{sec:proofs-related-s-1}

\begin{proofof}{Proposition~\ref{prop:regCont}}
  By~\ref{item:1notx} in Proposition~\ref{prop:1notx}, we may assume
  that
  $F_\nu (t,x,u) = \int_0^u \partial_u f (t,x,w) \, E'_\nu (w) \d{w}$
  and $F (t,x,u) = \int_0^u \partial_u f (t,x,w) \, E' (w)
  \d{w}$. Then, for all $(t,x) \in \spt\phi$ and all $v \in J_u$,
  \begin{eqnarray*}
    \modulo{F_\nu (t,x,v) - F (t,x,v)}
    & \leq
    & \modulo{
      \int_0^v \modulo{\partial_u f (t,x,w)} \,
      \modulo{E'_\nu (w) - E' (w)} \d{w}
      }
    \\
    & \leq
    & \mathcal{L} (J_u) \,
      \norma{\partial_u f}_{\L\infty (\spt\phi\times J_u;\reali)}
      \norma{E'_\nu - E'}_{\L\infty (J_u;\reali)} \,;
    \\
    \modulo{\div_x F_\nu (t,x,v) - \div_x F (t,x,v)}
    & \leq
    & \modulo{
      \int_0^v \modulo{\div_x \partial_u f (t,x,w)} \,
      \modulo{E'_\nu (w) - E' (w)} \d{w}
      }
    \\
    & \leq
    & \mathcal{L} (J_u) \,
      \norma{\div_x \partial_u f}_{\L\infty (\spt\phi\times J_u;\reali)}
      \norma{E'_\nu - E'}_{\L\infty (J_u;\reali)} \,,
  \end{eqnarray*}
  The proof follows by~\eqref{eq:37} and the Dominated Convergence Theorem.
\end{proofof}

\begin{proofof}{Theorem~\ref{thm:formulaC1}}
  We distinguish a few steps.

  \paragraph{Step~1. Equality~\eqref{eq:3} holds when
    $E_c (k) = \modulo{k-c}$, for a $c \in \reali$.}

  By~\eqref{eq:20}, Definition~\ref{def:muMnotx} and
  Definition~\ref{def:muk}, the left hand side in~\eqref{eq:3}
  satisfies
  \begin{displaymath}
    \left(\overline{\mathcal{M}}_{u} (E_c)\right) (\phi) = \mu_c (\phi) \,.
  \end{displaymath}
  To compute the right hand side in~\eqref{eq:3}, prepare $E_c'$
  following~\eqref{eq:20}:
  \begin{displaymath}
    E_c' (k) = \left\{
      \begin{array}{l@{\qquad}r@{\,}c@{\,}l}
        1
        & k
        & >
        & c
        \\
        0
        & k
        & =
        & c
        \\
        -1
        & k
        & <
        & c
      \end{array}
    \right.
    \quad\mbox{ equivalently } \quad
    \begin{array}{r@{\,}c@{\,}l@{\qquad}r@{\,}c@{\,}l}
      E_c'
      & =
      & 0 \vphantom{\caratt{[c,b]}}
      & c
      & \in
      & \reali \setminus [a,b]
      \\
      E_c'
      & =
      & -\caratt{[a,c]} + \caratt{[c,b]}
      & c
      & \in
      & [a,b]
      \\
      E_c'
      & =
      & \caratt{\mathopen]a,b]}
      & c
      & =
      &a
      \\
      E_c'
      & =
      & \caratt{[a,b\mathclose[}
      &
        c
      & =
      & b \,.
    \end{array}
  \end{displaymath}
  To complete this Step, apply~\cite[Lemma~6.3.2]{MR3839599}, also
  recalling Remark~\ref{rem:support}.

  \paragraph{Step~2: Extension to a general $E$.} Since both sides
  in~\eqref{eq:3} are linear in $E$, the result of Step~1 can be
  extended to any $E$ such that $E'$ is piecewise constant and
  satisfies~\eqref{eq:20}.

  By~\cite[(iii) in Exercise~4.2.5]{MR3839599}, we know that any
  $E \in \AC (\reali;\reali)$ with $E'$ regulated and
  satisfying~\eqref{eq:20} can be approximated by a sequence of
  functions $E_n$ such that $E_n'$ is piecewise constant and
  satisfies~\eqref{eq:20}, in the sense that both $E_n \to E$ and
  $E'_n \to E'$ uniformly on $J_u$.

  By~\cite[(i) in Theorem~6.3.9]{MR3839599} and~\ref{item:11} in
  Proposition~\ref{prop:2}, the right hand side in~\eqref{eq:3} is
  continuous with respect to the uniform convergence of $E$ and
  $E'$. The left hand side is continuous with respect to the uniform
  convergence of $E$ and $E'$ by Proposition~\ref{prop:regCont}. The
  proof is completed.
\end{proofof}

\begin{proofof}{Corollary~\ref{cor:solution}}
  The implication~\ref{item:24} $\implies$~\ref{item:26} follows
  from~\eqref{eq:3} in Theorem~\ref{thm:formulaC1}.  The
  implication~\ref{item:26} $\implies$~\ref{item:24} follows from
  $\mu_k = \overline{\mathcal{M}}_u (E_k)$ whenever
  $E_k (u) = \modulo{u-k}$, by Definition~\ref{def:muk}
  and~\eqref{eq:3}. The implication~\ref{item:24}
  $\implies$~\ref{item:25} follows from~\ref{item:22} in
  Theorem~\ref{thm:formulaC2}.

  Finally, we prove that~\ref{item:25} $\implies$~\ref{item:24}. Fix a
  positive test function
  $\phi \in \mathcal{D}(I \times \reali^n; \reali)$. By~\ref{item:11}
  in Proposition~\ref{prop:2}, the map $k \mapsto \mu_k (\phi)$ is of
  class $\BV$. Fix a point $\overline{k}$ of continuity of this map.
  Introduce a sequence of convex functions
  $E_\nu \in \C\infty (\reali;\reali)$ such that
  $E_\nu'' \underset{\nu \to+\infty}{\wsto} 2\,\delta_{\overline{k}}$,
  $\delta_{\overline{k}}$ being the Dirac $\delta$ centered at
  $\overline{k}$. By~\ref{item:25},
  $\left(\mathcal{M}_u (E_\nu)\right) (\phi) \geq 0$. By~\eqref{eq:2},
  with the same procedure as in~\cite[Theorem~4.1]{EvansGariepy2015},
  \begin{displaymath}
    0 \leq
    \left(\mathcal{M}_u (E_\nu)\right) (\phi)
    =
    \dfrac12 \int_{J_u} E_\nu'' (k) \; \mu_k (\phi) \; \d{k}
    \underset{\nu \to+\infty}{\longrightarrow}
    \mu_{\overline{k}} (\phi)
  \end{displaymath}
  by the above choices
  $E_\nu'' \underset{\nu \to+\infty}{\wsto} 2\,\delta_{\overline{k}}$
  and $\overline{k}$ point of continuity. Hence,
  $\mu_{\overline{k}} (\phi)\geq 0$. The case of $\overline{k}$ being
  a discontinuity point of $k \mapsto \mu_k (\phi)$ now follows from
  the previous case, thanks to~\ref{item:11} and~\ref{item:27}.
\end{proofof}

In a less general case, a proof of the
equivalence~\ref{item:24}~$\iff$~\ref{item:25} stated in
Corollary~\ref{cor:solution} is detailed
in~\cite[Proposition~2.4]{ColomboPerrollazDepX}.

\subsection{Proofs Related to \S~\ref{subsec:fourier-transform}}

\begin{proofof}{Proposition~\ref{prop:Fourier}}
  Consider the two items separately.

  \paragraph{\ref{item:19}:} Fix an arbitrary
  $\phi \in \mathcal{D} (I \times \reali^n; \complessi)$ and a
  $\xi \in \complessi$. Then, setting
  $U \colonequals \norma{u}_{\L\infty (I\times\reali^n;
    \reali)}$,
  \begin{flalign*}
    \left(\widehat\mu (\xi)\right) (\phi) %
    = %
    & \int_{\reali} \mu_k (\phi) \, e^{-i\, \xi \, k} \d{k}
    &[\mbox{By~\eqref{eq:15}}]
    \\
    = %
    & \int_{-U}^{U} \; \sum_{\nu=0}^{+\infty} \mu_k (\phi) \,
    \dfrac{(-i\, \xi k)^\nu}{\nu!}  \d{k}
    \\
    = %
    & \sum_{\nu=0}^{+\infty} \int_{-U}^{U} \mu_k (\phi) \,
    \dfrac{(-i\, \xi k)^\nu}{\nu!}  \d{k} %
    & [\mbox{By Uniform Convergence}]
    \\
    = %
    & \sum_{\nu=0}^{+\infty} \dfrac{(-i)^\nu}{\nu!}  \int_{-U}^{U}
    k^\nu \, \mu_k (\phi) \d{k} \; \xi^\nu \,.
  \end{flalign*}
  Thanks to~\ref{item:11} in Proposition~\ref{prop:2}, this shows that
  $\xi \mapsto \left(\widehat\mu (\xi)\right) (\phi)$ is expressed as
  a power series with infinite radius of convergence.

  \paragraph{\ref{item:20}:} By construction,
  $\mathcal{E}_\xi \in \C\infty (\reali; \complessi)$ and
  $\mathcal{E}_\xi'' (k) = 2\,e^{-i\, \xi \,k}$. The proof follows by an
  application of~\eqref{eq:2} in Theorem~\ref{thm:formulaC2}.
\end{proofof}

\begin{proofof}{Theorem~\ref{thm:propFT}}
  Introduce $\mu_k$ as in Definition~\ref{def:muk} and call $\nu_k$
  the analogous distribution referred to $v$ and
  $v_o$. By~\ref{item:22} in Theorem~\ref{thm:formulaC2},
  $\mathcal{M}_{u,u_o} \left(\mathcal{T}_kE \right) =
  \mathcal{M}_{v,v_o} \left(\mathcal{T}_kE\right)$ implies that for
  all $k \in \reali$ and all
  $\phi \in \mathcal{D} (I \times \reali^n; \complessi)$
  \begin{displaymath}
    \int_{\reali} E'' (w-k) \; \mu_w (\phi) \d{w}
    =
    \int_{\reali} E'' (w-k) \; \nu_w (\phi) \d{w}
  \end{displaymath}
  also using the fact that $\mu_k$ and $\nu_k$ are compactly
  supported. Introduce $q (p) \colonequals E'' (-p)$ and
  $\Delta (w) \colonequals (\mu_w-\nu_w) (\phi)$ so that the latter
  equality can be rewritten as
  \begin{displaymath}
    q * \Delta = 0 \,.
  \end{displaymath}
  The left hand side above is a convolution between the continuous
  real function $q$ and the compactly supported $\BV$ function
  $\Delta$ attaining values in $\complessi$. Since $q$ is a
  tempered distribution and $\Delta$ is a compact distribution by
  Remark~\ref{rem:support} and~\cite[Theorem~7.1.15]{MR1065993} we get
  \begin{displaymath}
    \widehat{q} \; \widehat{\Delta} =
    \widehat{(q * \Delta)} =
    0 \,.
  \end{displaymath}
  Here, $\widehat{q}$ is a tempered distribution and $\widehat\Delta$
  is an entire function. Thus, for any
  $\psi \in \mathcal{D} (\reali;\complessi)$
  \begin{displaymath}
    \widehat{q} (\widehat{\Delta} \; \psi) = 0
  \end{displaymath}
  Hence, $\widehat\Delta$ vanishes on all $\spt \widehat{q}$
  (understood in the sense
  of~\cite[Definition~2.2.2]{MR1065993}). Since $\spt \widehat {E''}$
  admits a cluster point, the set of zeros of $\widehat\Delta$ admits
  a cluster point.  Moreover, $\widehat{\Delta}$ is an entire function
  by~\ref{item:19} in Proposition~\ref{prop:Fourier}, so it vanishes
  either everywhere or only in isolated zeros. Thus,
  $\widehat\Delta = 0$. Hence also $\Delta = 0$. By the arbitrariness
  of $\phi$, for all $w$ we have $\mu_w = \nu_w$. The proof is
  completed thanks to~\ref{item:22} in Theorem~\ref{thm:formulaC2}.
\end{proofof}

\appendix

\section{Integrals in \texorpdfstring{$\mathcal{D}'$}{D'}}
\label{sec:integr-texorpdfstr}

For completeness, we state below the following lemma on integrals
taking values in a distribution space, refer to~\cite{zbMATH03576139}
for related results.

\begin{proposition}
  \label{prop:integral}
  Let $\Omega \subseteq \reali^d$ be non empty and either open or
  $\overline{\Omega} = \overline{\pint{\Omega}}$ is the closure of an
  open set with $\C1$ boundary. Fix a finite measure space
  $(X, \mathcal{S}, \mathfrak{m})$ and let
  $\Upsilon \colon X \to \mathcal{D}' (\Omega;\complessi)$.  If
  \begin{enumerate}[label=(\arabic*)]
  \item \label{item:12} For all
    $\phi \in \mathcal{D} (\Omega;\complessi)$, the map
    \begin{equation}
      \label{eq:7}
      \begin{array}{ccccc}
        \Upsilon (\phi)
        & \colon
        & X
        & \to
        & \complessi
        \\
        &
        & x
        & \mapsto
        & \left(\Upsilon (x)\right) (\phi)
      \end{array}
    \end{equation}
    is measurable.
  \item \label{item:13}For all compact subsets $K$ of $\Omega$, there
    exist a positive $C$ and an integer $n$ such that for all
    $\phi \in \mathcal{D} (\Omega;\complessi)$ with
    $\spt \phi \subseteq K$, we have that for all $x \in X$
    \begin{equation}
      \label{eq:8}
      \modulo{\left(\Upsilon (x)\right) (\phi)}
      \leq
      C \, \norma{\phi}_{\C{n} (K;\complessi)} \,.
    \end{equation}
  \end{enumerate}
  \noindent Then, the map
  \begin{equation}
    \label{eq:9}
    \begin{array}{ccccc}
      \int_X \Upsilon \d{\mathfrak{m}}
      & \colon
      & \mathcal{D} (\Omega;\complessi)
      & \to
      & \complessi
      \\
      &
      & \phi
      & \mapsto
      & \int_X \left(\Upsilon (x)\right) (\phi) \, \d{\mathfrak{m} (x)}
    \end{array}
  \end{equation}
  is well defined and belongs to $\mathcal{D}' (\Omega;\complessi)$.
\end{proposition}

\begin{proofof}{Proposition~\ref{prop:integral}}
  By~\eqref{eq:7} and~\eqref{eq:8}, since $X$ has finite measure, for
  any $\phi \in \mathcal{D} (\Omega;\complessi)$, the quantity
  $\int_X \left(\Upsilon (x)\right) (\phi) \, \d{\mathfrak{m} (x)}$ is
  a well defined Lebesgue integral. Hence, the map in~\eqref{eq:9} is
  well defined and linear.

  Fix a compact $K$ in $\Omega$. Then select $C$ and $n$ as
  in~\ref{item:12} and using $\mathfrak{m} (X) < +\infty$, obtain that
  for any $\phi \in \mathcal{D}' (\Omega;\complessi)$ with
  $\spt \phi \subseteq K$,
  $\modulo{\int_X \left(\Upsilon (x)\right) (\phi) \, \d{\mathfrak{m}
      (x)}} \leq C \, \mathfrak{m} (X) \, \norma{\phi}_{\C{n}
    (K;\complessi)}$, completing the proof.
\end{proofof}

In \S~\ref{subsec:case-class-entr} the following extension of the
Kurzweil-Stieltjes integral to distribution valued functions is
necessary.

\begin{proposition}
  \label{prop:integral2}
  Let $\Omega \subseteq \reali^d$ be non empty and open. Fix a real
  interval $[a,b]$ with $a<b$. Let $\Gamma \in \BV ([a,b];\reali)$ and
  $\Upsilon \colon [a,b] \to \mathcal{D}' (\Omega;\reali)$.  If
  \begin{enumerate}[label=(\arabic*)]
  \item \label{item:14} For all
    $\phi \in \mathcal{D} (\Omega;\reali)$, the map
    \begin{equation}
      \label{eq:10}
      \begin{array}{ccccc}
        \langle \Upsilon | \phi \rangle
        & \colon
        & [a,b]
        & \to
        & \reali
        \\
        &
        & x
        & \mapsto
        & \left(\Upsilon (x)\right) (\phi)
      \end{array}
    \end{equation}
    is in $\BV ([a,b];\reali)$.
  \item \label{item:15} For all compact subsets $K$ of $\Omega$, there exist a
    positive $C$ and an integer $n$ such that for all
    $\phi \in \mathcal{D} (\Omega;\reali)$ with
    $\spt \phi \subseteq K$, we have that for all $x \in [a,b]$
    \begin{equation}
      \label{eq:11}
      \modulo{\left(\Upsilon (x)\right) (\phi)}
      \leq
      C \, \norma{\phi}_{\C{n} (K;\reali)} \,.
    \end{equation}
  \end{enumerate}
  \noindent Then, the map
  \begin{equation}
    \label{eq:9bis}
    \begin{array}{ccccc}
      \int_a^b \Upsilon \d{\Gamma}
      & \colon
      & \mathcal{D} (\Omega;\reali)
      & \to
      & \reali
      \\
      &
      & \phi
      & \mapsto
      & \int_a^b \left(\Upsilon (x)\right) (\phi) \, \d{\Gamma (x)}
    \end{array}
  \end{equation}
  is well defined and belongs to $\mathcal{D}' (\Omega;\reali)$.
\end{proposition}

\begin{proofof}{Proposition~\ref{prop:integral2}}
  The integral in the right hand side of~\eqref{eq:9bis} is understood in
  the Kurzweil-Stieltjes sense, see~\cite[Chapter~6]{MR3839599} and in
  particular~\cite[Definition~6.2.2]{MR3839599}. Note that, for any
  $\phi \in \mathcal{D} (\Omega;\reali)$, the same integral is well
  defined by~\cite[Theorem~6.3.11]{MR3839599}, which can be applied
  by~\ref{item:14} and since $\Gamma \in \BV ([a,b];\reali)$. This
  ensures that the map $\int_a^b \Upsilon \d{\Gamma}$ defined
  in~\eqref{eq:9bis} is well defined. Moreover, it is linear
  by~\cite[Theorem~6.2.7]{MR3839599}.

  Fix a compact $K \subseteq \Omega$. Then select $C$ and $n$ as
  in~\ref{item:14} and obtain by~\cite[Theorem~6.3.6]{MR3839599} that
  for any $\phi \in \mathcal{D}' (\Omega;\reali)$ with
  $\spt \phi \subseteq K$,
  $\modulo{\int_X \left(\Upsilon (x)\right) (\phi) \, \d{\Gamma (x)}}
  \leq C \; \tv (\Gamma) \, \norma{\phi}_{\C{n} (K;\reali)}$,
  completing the proof.
\end{proofof}

\section{Tensor Products of \texorpdfstring{$\C1$}{C1} Spaces}
\label{sec:tensor-product-c1}

We collect below some elementary properties of $\mathscr{F}$,
$\mathscr{F}_c$ and of their topologies, as defined
in~\eqref{eq:42}, \eqref{eq:45}, \eqref{eq:41} and use the notation~\eqref{eq:50}. For the basic definitions and
terminology, we refer to~\cite[Chapters~7, 10 and~13]{MR2296978}.

\begin{proposition}
  \label{prop:tensor}
  The space $\mathscr{F}$, with the seminorms
  $\left\{p_K\colon K \subset I \times \reali^n \times \reali \mbox{
      compact}\right\}$ is a locally convex Hausdorff vector space and
  a Fr\'echet space. Moreover, if $\zeta\in \mathscr{F}$, then also the
  mixed derivatives $\partial^2_{ut}\zeta$ and $\partial_u \grad_x\zeta$
  exist, are continuous on $I \times \reali^n \times \reali$ and
  \begin{equation}
    \label{eq:43}
    \partial^2_{ut}\zeta = \partial^2_{tu}\zeta
    \qquad\qquad
    \partial_u \grad_x\zeta = \grad_x \partial_u \zeta \,.
  \end{equation}
\end{proposition}

\begin{proofof}{Proposition~\ref{prop:tensor}}
  It is immediate to prove that $\mathscr{F}$ is a locally convex
  Hausdorff metrizable vector space. Note that any function $\zeta$ in
  $\mathscr{F}$ satisfies for all
  $(t,x,u), (t',x',u') \in I \times \reali^n \times \reali$
  \begin{eqnarray*}
    \zeta (t',x,u')
    & =
    & \zeta (t,x,u)
      + \int_t^{t'} \partial_t \zeta (s,x,u) \d{s}
      + \int_u^{u'} \partial_u \zeta (t,x,w) \d{w}
    \\
    &
    & \quad
      +\int_u^{u'} \int_t^{t'} \partial^2_{tu} \zeta (s,x,w) \d{s} \d{w} \,;
    \\
    \zeta (t,x+ (x_i'-x_i)\mathbf{e}_i,u')
    & =
    & \zeta (t,x,u)
      + \int_u^{u'} \partial_u \zeta (t,x,w) \d{w}
      + \int_{0}^{x_i'-x_i} \partial_{x_i} \zeta (t,x+y\,\mathbf{e}_i,u) \d{y}
    \\
    &
    & \quad
      +\int_u^{u'} \int_{0}^{x_i'-x_i} \partial^2_{x_iu} \zeta (t,x+y\,\mathbf{e}_i,w) \d{y} \d{w}
  \end{eqnarray*}
  for $i=1, \ldots, n$, $\mathbf{e}_1, \ldots, \mathbf{e}_n$ being the
  canonical base in $\reali^n$. The above relations directly imply the
  completeness of $\mathscr{F}$ with respect to the topology induced
  by the seminorms. The existence, continuity of $\partial^2_{ut}\zeta$
  and $\partial_u \grad_x \zeta$ as well as~\eqref{eq:43} also follow,
  thanks to Fubini Theorem and to the Fundamental Theorem of Calculus.
\end{proofof}

\begin{lemma}
  \label{lem:treves}
  $\mathscr{F}_c$ is an inductive limit of Fr\'echet spaces.  A
  linear map on $\mathscr{F}_c$ is continuous if and only if it is
  sequentially continuous.
  $\zeta_\nu \underset{\nu\to+\infty}{\longrightarrow} \zeta$ in
  $\mathscr{F}_c$ if and only if
  $\zeta_\nu \underset{\nu\to+\infty}{\longrightarrow} \zeta$ in
  $\mathscr{F}$ and there exists a compact
  $K_{tx} \subset I \times \reali^n$ such that
  ${\spt}_{tx} \zeta_\nu \subseteq K_{tx}$ for all $\nu$.
\end{lemma}

\begin{proofof}{Lemma~\ref{lem:treves}}
  Let $r_\nu$ be an increasing sequence converging to $\sup I$. Note
  that
  \begin{displaymath}
    \mathscr{F}_c
    =
    \bigcup_{\nu \in \naturali} \mathscr{F}_\nu
    \quad \mbox{ where } \quad
    \mathscr{F}_\nu
    \colonequals
    \left\{
      \zeta \in \mathscr{F} \colon
      {\spt}_{tx} \zeta \subset  [0,r_\nu] \times \overline{B_{\reali^n} (0,\nu)}
    \right\}
  \end{displaymath}
  and $\mathscr{F}_\nu$ is a Fr\'echet space, since it is a closed
  subset of $\mathscr{F}$, by Proposition~\ref{prop:tensor}. Hence,
  by~\cite[Proposition~13.1]{MR2296978}, the continuity of a linear
  map on $\mathscr{F}_c$ is equivalent to the continuity of its
  restrictions to each $\mathscr{F}_\nu$. The latter statement follows
  from~\cite[Corollary~1, Chapter~14]{MR2296978}.
\end{proofof}

We now show that the algebraic tensor product
$\C1 (I \times \reali^n; \reali) \otimes \C1 (\reali;\reali)$ is
dense in $\mathcal{F}$ with respect to the topology induced
by~\eqref{eq:45}. This is achieved by a refinement of the classical
{Weierstra\ss} Theorem.

\begin{proposition}
  \label{prop:density}
  Fix $\zeta \in \mathscr{F}$. For any compact set
  $K \subset I \times \reali^n \times \reali$ and for any
  $\epsilon > 0$, there exists $N \in \naturali$,
  $\phi_1, \ldots, \phi_N \in \C\infty (I \times \reali^n; \reali)$ and
  $e_1, \ldots, e_N \in \C\infty (\reali;\reali)$ such that
  \begin{equation}
    \label{eq:46}
    p_K\left(\zeta - \sum_{i=1}^N \phi_i \otimes e_i \right) < \epsilon \,.
  \end{equation}
\end{proposition}

Above, we used the standard tensor product notation: $(\phi_i \otimes e_i) (t,x,u) = \phi_i (t,x) \; e_i (u)$.

\begin{proofof}{Proposition~\ref{prop:density}}
  It is sufficient to consider the case
  $K = [-R,R] \times [-R,R]^n \times [-R,R]$, for a positive $R$.

  Introduce a function
  $\tilde\zeta \in \C1 (\reali\times\reali^n\times\reali;\reali)$ such
  that
  $\partial^2_{tu} \tilde\zeta \in
  \C0(\reali\times\reali^n\times\reali;\reali)$,
  $\grad_x \partial_u \tilde\zeta \in
  \C0(\reali\times\reali^n\times\reali;\reali^n)$ and
  $\rest{\tilde\zeta}{I\times\reali^n\times\reali} = \zeta$.

  In the case $I=\reali_+$, a possible construction of $\tilde\zeta$ is
  \begin{displaymath}
    \tilde\zeta (t,x,u)
    \colonequals
    \left\{
      \begin{array}{l@{\qquad}r@{\,}c@{\,}l}
        \zeta (t,x,u)
        & (t,x,u)
        & \in
        & I \times \reali^n \times \reali
        \\
        2\zeta (0,x,u) - \zeta (-t,x,u)
        & (t,x,u)
        & \in
        & \mathopen]-\infty,0\mathclose[ \times \reali^n \times \reali \,.
      \end{array}
      \right.
  \end{displaymath}
  The same reflection techniques, coupled with an induction procedure,
  allows to tackle the case $I$ is bounded.

  Then, multiply $\tilde\zeta$ by a suitable smooth function to obtain
  a map
  $\overline{\zeta} \in \C1 (\reali \times \reali^n \times \reali;
  \reali)$ such that
  $\partial^2_{tu} \overline{\zeta} \in
  \C0(\reali\times\reali^n\times\reali;\reali)$,
  $\grad_x \partial_u \overline{\zeta} \in
  \C0(\reali\times\reali^n\times\reali;\reali^n)$ with moreover
  \begin{displaymath}
    \rest{\overline{\zeta}}{K} = \rest{\tilde\zeta}{K}
    \quad \mbox{ and } \quad
    \spt \overline{\zeta} \subseteq 2K \,.
  \end{displaymath}
  Define the regularizing kernel
  \begin{equation}
    \label{eq:49}
    \!\!\!
    \rho_\nu (t,x,u)
    \colonequals
    \overline\rho_\nu (t) \; \overline\rho_\nu (u)
    \prod_{i=1}^n \overline\rho_\nu (x_i)
    \quad \mbox{ where } \quad
    \overline\rho_\nu (\xi) \colonequals
    \left\{
      \begin{array}{@{}lr@{\,}c@{\,}l}
        c_\nu \, \left(1-\left(\frac{\xi}{3R}\right)^2\right)^\nu
        & \modulo{\xi}
        & \leq
        & 3R
        \\
        0
        & \modulo{\xi}
        & >
        & 3R
      \end{array}
    \right.
  \end{equation}
  and $c_\nu$ is chosen so that
  $\int_{\reali} \overline\rho_\nu (\xi)\d\xi = 1$. Let
  $\zeta_\nu \colonequals \overline{\zeta} * \rho_\nu$. For all
  $(t,x,u) \in K$ we have
  \begin{eqnarray*}
    &
    & \zeta_\nu (t,x,u)
    \\
    & =
    & \int_{2K} \overline{\zeta} (s,y,w) \;
      \rho_n (t-s, x-y, u-w) \d{s} \d{y} \d{w}
    \\
    & =
    & \! {c_\nu}^{n+2} \!
      \int_{2K} \! \overline{\zeta} (s,y,w)
      \left[
      \left(1{-}\left(\frac{t{-}s}{3R}\right)^2\right) \!
      \left(1{-}\left(\frac{u{-}w}{3R}\right)^2\right) \!
      \prod_{i=1}^n
      \left(1{-}\left(\frac{x_i{-}y_i}{3R}\right)^2\right) \!
      \right]^\nu
      \d{s} \d{y} \d{w} \,.
  \end{eqnarray*}
  The latter equality holds because $(t,x,u) \in K$, $(s,y,w) \in 2K$
  so that $(t-s, x-y,u-w) \in 3K$ and~\eqref{eq:49} applies. Hence,
  $\rest{\zeta_\nu}K$ coincides with a polynomial in $(t,x,u)$ that
  can be written in the form $\sum_{i=1}^N \phi_i \otimes e_i$, where
  all the $\phi_i$ and $e_i$ are of class $\C\infty$.

  By the standard properties of the convolution,
  see~\cite[\S~1.3]{MR1065993},
  \begin{displaymath}
    \begin{array}{rcl}
      \partial_t \zeta_\nu
      & =
      & (\partial_t \overline{\zeta}) * \rho_\nu
      \\
      \grad_x \zeta_\nu
      & =
      & (\grad_x \overline{\zeta}) * \rho_\nu
      \\
      \partial_u \zeta_\nu
      & =
      & (\partial_u \overline{\zeta}) * \rho_\nu
    \end{array}
    \quad \mbox{ and } \quad
    \begin{array}{rcl}
      \partial^2_{tu} \zeta_\nu
      & =
      & (\partial^2_{tu} \overline{\zeta}) * \rho_\nu
      \\
      \grad_x \partial_u \zeta_\nu
      & =
      & (\grad_x \partial_u \overline{\zeta}) * \rho_\nu
    \end{array}
  \end{displaymath}
  so that $\lim_{\nu\to +\infty} p_K (f-\zeta_\nu) = 0$.
\end{proofof}

Remark that the above result actually provides an elementary
description of the projective tensor product of
$\C1 (I \times \reali^n; \reali) \otimes \C1 (\reali;\reali)$. For
the general theory, refer to~\cite[Section~45, Theorem~45.1
and Corollary~3]{MR2296978}.

\begin{corollary}
  \label{cor:density}
  The algebraic tensor product
  $\Cc\infty (I \times \reali^n; \reali) \otimes \C\infty
  (\reali;\reali)$ is dense in $\mathcal{F}_c$ with respect to the
  inductive topology.
\end{corollary}

\begin{proofof}{Corollary~\ref{cor:density}}
  Fix $\zeta \in \mathscr{F}_c$. Since
  $\mathscr{F}_c \subset \mathscr{F}$, by
  Proposition~\ref{prop:density}, for all $\nu \in \naturali$ there
  exist
  $\zeta_\nu \in \C\infty (I\times\reali^n;\reali) \otimes \C\infty
  (\reali;\reali)$ such that
  $\zeta_\nu \underset{\nu\to+\infty}{\longrightarrow} \zeta$ in
  $\mathscr{F}$. Introduce a map
  $\phi \in \Cc\infty (I \times \reali^n;\reali)$ such that
  $\phi (t,x) = 1$ for all $t,x \in {\spt}_{tx}\zeta$. Then, in the
  topology of $\mathscr{F}_c$,
  $\phi\,\zeta_\nu \underset{\nu\to+\infty}{\longrightarrow} \phi\,
  \zeta = \zeta$, completing the proof.
\end{proofof}

\begin{proposition}
  \label{prop:fEF}
  Consider the map
  \begin{equation}
    \label{eq:44}
    \begin{array}[b]{ccccc}
      \mathcal{F}
      &\colon
      & \mathscr{F}_c \times \mathscr{F}^n
      & \to
      & \mathscr{F}_c^n
      \\
      &
      & (E,f)
      & \mapsto
      & F
    \end{array}
    \mbox{ where }\;
    F (t,x,u)
    \colonequals
    \int_0^u \partial_u E (t,x,w) \; \partial_u f(t,x,w) \d{w}
  \end{equation}
  Then,
  \begin{enumerate}[label=($\mathcal{F}\arabic*$)]

  \item \label{item:38} $\mathcal{F}$ is well defined.
  \item \label{item:39} $\mathcal{F}$ is bilinear.
  \item \label{item:43} $\mathcal{F}$ is continuous, in the product topology.
  \end{enumerate}
\end{proposition}

\noindent The introduction of the map $\mathcal{F}$ is motivated by
Definition~\ref{def:EEF}: $\mathcal{F} (E,f)$ is an entropy flux
paired to the entropy $E$ with respect to the flux $f$.

\begin{proofof}{Proposition~\ref{prop:fEF}}
  Note first that
  ${\spt}_{tx} \mathcal{F} (E,f) \subseteq {\spt}_{tx} E$.  By direct
  computations,
  \begin{eqnarray*}
    \partial_t F (t,x,u)
    & =
    & \int_0^u \partial^2_{tu} E (t,x,w) \; \partial_u f(t,x,w) \d{w}
      +
      \int_0^u \partial_u E (t,x,w) \; \partial^2_{tu} f(t,x,w) \d{w}
    \\
    \grad_x F (t,x,u)
    & =
    & \int_0^u \grad_x \partial_u E (t,x,w) \; \partial_u f(t,x,w) \d{w}
    \\
    &
    & \qquad
      +
      \int_0^u \partial_u E (t,x,w) \; \grad_x \partial_u f(t,x,w) \d{w}
    \\
    \partial_u F (t,x,u)
    & =
    & \partial_u E (t,x,w) \; \partial_u f(t,x,w)
    \\
    \partial^2_{tu} F (t,x,u)
    & =
    & \partial^2_{tu} E (t,x,w) \; \partial_u f(t,x,w)
      +
      \partial_u E (t,x,w) \; \partial^2_{tu} f(t,x,w)
    \\
    \grad_x \partial_u F (t,x,u)
    & =
    & \grad_x \partial_u E (t,x,w) \; \partial_u f(t,x,w)
      +
      \partial_u E (t,x,w) \; \grad_x \partial_u f(t,x,w)
  \end{eqnarray*}
  showing that the range of $\mathcal{F}$ is in $\mathscr{F}_c$. The
  bilinearity~\ref{item:39} is immediate.

  To prove the continuity~\ref{item:43} of $\mathcal{F}$, by
  Lemma~\ref{lem:treves} it is sufficient to prove sequential
  continuity. Let $E_\nu$, $f_\nu$ be such that
  $E_\nu \underset{\nu\to+\infty}{\longrightarrow} E$ in
  $\mathscr{F}_c$ and
  $f_\nu \underset{\nu\to+\infty}{\longrightarrow} f$ in
  $\mathscr{F}$. Then,
  ${\spt}_{tx}\mathcal{F} (E_\nu,f_\nu) \subseteq {\spt}_{tx} E_\nu$
  and
  $p_K \left(\mathcal{F} (E_\nu,f_\nu) - \mathcal{F} (E,f)\right)
  \underset{\nu\to+\infty}{\longrightarrow} 0$ thanks to the
  {formul\ae} above. By Lemma~\ref{lem:treves}, the proof
  of~\ref{item:43} is completed.
\end{proofof}

\medskip

\noindent\textbf{Acknowledgment:} The first
author acknowledges the PRIN 2022 project \emph{Modeling, Control and
  Games through Partial Differential Equations} (D53D23005620006),
funded by the European Union - Next Generation EU and the 2024~GNAMPA
project \emph{Modelling and Analysis through Conservation Laws}. The
second author was supported by the ANR project COSS \emph{COntrol on
  Stratified Structures} (ANR-22-CE40-0010) and from the \emph{Fund
  for International Cooperation} from the University of Brescia.

\bigskip

\paragraph{Data Availability Statement.} No data sets were generated
or analyzed during the current study.

\paragraph{Conflict of Interest Statement.} The authors have no
competing interests to declare that are relevant to the content of
this article.

{

  \small

  \bibliography{colombo_perrollaz-entropy_production}

  \bibliographystyle{abbrv}

}

\end{document}